\documentclass[10pt,sumlimits,namelimits]{article}

\usepackage{amsthm}
\usepackage{amsmath}
\usepackage{amsfonts}
\usepackage{amssymb}
\usepackage{mathrsfs}
\usepackage{stmaryrd}
\usepackage{dsfont}
\usepackage{enumerate}
\usepackage{color}
\usepackage{hyperref}
\usepackage[applemac]{inputenc}

\newcommand{\1}{\mathds{1}}
\newcommand{\A}{\mathbb{A}}

\newcommand{\N}{\mathbb{N}}

\newcommand{\R}{\mathbb{R}}
\newcommand{\C}{\mathbb{C}}
\renewcommand{\O}{\Omega}
\renewcommand{\o}{\omega}

\newcommand{\vO}{{\cal O}}

\newcommand{\vR}{{\cal R}}

\newcommand{\Dr}{\mathscr{D}}

\newcommand{\vphi}{\varphi}
\newcommand{\eps}{\varepsilon}
\newcommand{\la}{\lambda}
\newcommand{\dsp}{\displaystyle}
\newcommand{\ovl}{\overline}
\newcommand{\udl}{\underline}

\newcommand{\vsup}{\sup\limits}

\newcommand{\vint}{\int\limits}
\newcommand{\vsum}{\sum\limits}
\newcommand{\inj}{\hookrightarrow}
\newcommand{\tends}{\longrightarrow}
\newcommand{\weak}{\rightharpoonup}
\newcommand{\wt}{\widetilde}

\newcommand{\loc}{\mathrm{loc}}
\newcommand{\rad}{\mathrm{rad}}
\newcommand{\CP}{C_\mathrm{P}}

\renewcommand{\c}{\mathrm{c}}
\renewcommand{\d}{\mathrm{d}}

\newcommand{\dist}{\mathrm{dist}}

\newcommand{\w}{{\textsl w}}
\renewcommand{\le}{\leqslant}
\renewcommand{\ge}{\geqslant}
\renewcommand{\Re}{\mathrm{Re}}
\renewcommand{\Im}{\mathrm{Im}}
\newcommand{\bs}{\boldsymbol}
\newcommand{\vi}{\mathrm{i}}
\newcommand{\p}{\prime}
\newcommand{\vecteur}{\overrightarrow}
\DeclareMathOperator{\supp}{supp}

\numberwithin{equation}{section}

\newtheorem{thm}{Theorem}[section]
\newtheorem{prop}[thm]{Proposition}
\newtheorem{cor}[thm]{Corollary}
\newtheorem{lem}[thm]{Lemma}
\newtheorem{sym}[thm]{Symmetry Property}

\theoremstyle{definition}
\newtheorem{rmk}[thm]{Remark}
\newtheorem{defi}[thm]{Definition}
\newtheorem{com}[thm]{Comments}

\newenvironment{proof*}{\noindent{\bf Proof.}}{\qed}
\newenvironment{vproof}[1]{\noindent{\bf Proof #1}}{\qed}

\title{\huge \sc Existence of weak solutions to some stationary Schrödinger equations with singular nonlinearity}
\author{\sc Pascal Bégout$^*$ and Jes\'us Ildefonso D\'iaz$^\dagger$}
\date{}

\begin{document}

\maketitle

\begin{gather*}
\begin{array}{cc}
            ^*\mbox{Institut de Mathématiques de Toulouse \& TSE }	&	\;^\dagger\mbox{Instituto de Matem\'atica Interdisciplinar}		\\
                                          \mbox{Université Toulouse I Capitole }	&	\mbox{ Departamento de Matem\'atica Aplicada}			\\
                                                   \mbox{Manufacture des Tabacs }	&	\mbox{ Universidad Complutense de Madrid}				\\
                                                          \mbox{21, Allée de Brienne }	&	\mbox{ Plaza de las Ciencias, 3}						\\
                                  \mbox{31015 Toulouse Cedex 6, FRANCE }	&	\mbox{ 28040 Madrid, SPAIN}
\bigskip \\
\mbox{
{\footnotesize $^*$e-mail\:: \href{mailto:Pascal.Begout@math.cnrs.fr}{\udl{\texttt{Pascal.Begout@math.cnrs.fr}}}}}
&
\mbox{
{\footnotesize $^\dagger$e-mail\:: \href{mailto:diaz.racefyn@insde.es}{\udl{\texttt{diaz.racefyn@insde.es}}}}
}
\end{array}
\end{gather*}

\begin{abstract}
We prove some existence (and sometimes also uniqueness) of solutions to some stationary equations associated to the complex Schrödinger operator under the presence of a singular nonlinear term. Among other new facts, with respect some previous results in the literature for such type of nonlinear potential terms, we include the case in which the spatial domain is possibly unbounded (something which is connected with some previous localization results by the authors), the presence of possible non-local terms at the equation,  the case of boundary conditions different to the Dirichlet ones and, finally, the proof of the existence of solutions when the right-hand side term of the equation is beyond the usual $L^2$-space.
\end{abstract}

{\let\thefootnote\relax\footnotetext{$^\dagger$The research of J.I.~D\'iaz was partially supported by the project ref. MTM2011-26119 of the DGISPI (Spain) and the Research Group MOMAT (Ref. 910480) supported by UCM. He has received also support from the ITN \textit{FIRST} of the Seventh Framework Program of the European Community's (grant agreement number 238702)}}
{\let\thefootnote\relax\footnotetext{2010 Mathematics Subject Classification: 35Q55, (35A01, 35A02, 35B45, 35B65, 35J60)}}
{\let\thefootnote\relax\footnotetext{Key Words: nonlinear Schrödinger equation, different boundary conditions, unbounded domains, non local terms, data in weighted spaces, existence, uniqueness, smoothness}}

\tableofcontents

\baselineskip .65cm

\section{Introduction}
\label{intro}

\noindent
This paper is concerned by existence of solutions for two kinds of equations related to the complex Schrödinger operator,
\begin{gather}
\label{1}
-\Delta u+a|u|^{-(1-m)}u+bu=F, \text{ in } L^2(\O), \\
\label{2}
-\Delta u+a|u|^{-(1-m)}u+bu+cV^2u=F, \text{ in } L^2(\O),
\end{gather}
with homogeneous Dirichlet boundary condition
\begin{gather}
	\label{dir}
		 u_{|\Gamma}=0,
\end{gather}
or homogeneous Neumann boundary condition
\begin{gather}
	\label{neu}
		 \dfrac{\partial u}{\partial\nu}_{|\Gamma}=0,
\end{gather}
where $\O$ is a subset of $\R^N$ with boundary $\Gamma,$ $0<m<1,$ $(a,b,c)\in\C^3$ and $V\in L^\infty(\O;\R)$ is a real potential. Here and in what follows, when $\Gamma$ is of class $C^1,$ $\nu$ denotes the outward unit normal vector to $\Gamma.$ Moreover, $\Delta=\sum\limits_{j=1}^N\frac{\partial^2}{\partial x^2_j}$ is the Laplacian in $\O.$
\medskip \\
In Bégout and D\'iaz~\cite{MR2876246}, the authors study the spatial localization property compactness of the support of solutions of
equation~\eqref{1} (see Theorems~3.1, 3.5, 3.6, 4.1, 4.4 and 5.2). Existence, uniqueness and \textit{a priori} bound are also established with the homogeneous Dirichlet boundary condition, $F\in L^p(\O)$ $(2<p<\infty)$ and $(a,b)\in\C^2$ satisfying assumptions~\eqref{ab} below. In this paper, we give such existence and \textit{a priori} bound results but for the weaker assumption $F\in L^2(\O)$ (Theorems~\ref{thmexista3} and \ref{thmbound3}) and also for some different hypotheses on $(a,b)\in\C^2$ (Theorems~\ref{thmexista1} and \ref{thmbound1}). Additionally, we consider homogeneous Neumann boundary condition (Theorems~\ref{thmexista3} and \ref{thmbound3}).
\medskip \\
In Bégout and D\'iaz~\cite{MR3193996}, spatial localization property for the partial differential equation \eqref{2} associated to self-similar solutions of the nonlinear Schrödinger equation
\begin{gather*}
\vi u_t+\Delta u=a|u|^{-(1-m)}u+f(t,x),
\end{gather*}
is studied.
\medskip \\
In this paper, we prove existence of solutions with homogeneous Dirichlet or Neumann boundary conditions (Theorems~\ref{thmexista2}) and establish \textit{a priori} bounds (Theorem~\ref{thmbound2}), for both equations~\eqref{1} and \eqref{2} with any of both boundary conditions~\eqref{dir} or \eqref{neu}. We also show uniqueness (Theorem~\ref{thmuni}) and regularity results (Theorem~\ref{thmreg}), under suitable additional conditions. We send the reader to the long introduction of Bégout and D\'iaz~\cite{MR3193996} for many comments on the frameworks in which the equation arises (Quantum Mechanics, Nonlinear Optics and Hydrodynamics) and their connections with some other papers in the literature.
\medskip \\
This paper is organized as follows. In the next section, we give results about existence, uniqueness, regularity and \textit{a priori} bounds for
equations~\eqref{1} and \eqref{2}, with boundary conditions~\eqref{dir} or \eqref{neu}, and notations are given in Section~\ref{not}. Section~\ref{apriori}, is devoted to the establishment of \textit{a priori} bounds for the different truncated nonlinearities of equations studied in this paper. In
Section~\ref{proofs}, we prove the results given in Section~\ref{main}. In Bégout and D\'iaz~\cite{MR2876246}, localization property is studied for equation~\eqref{1}. The results we give require, sometimes, the same assumptions on $(a,b)\in\C^2$ as in Bégout and D\'iaz~\cite{MR2876246} but with a change of notation. See Comments~\ref{com} below for the motivation of this change. In Section~\ref{diaz} we will show the existence of solutions to equation~\eqref{2} for data in a weighted subspace. Finally, in the last section, we state the principal results obtained in this paper and give some applications. Existence of solutions for equation~\eqref{2} is used in Bégout and D\'iaz~\cite{MR3193996} while existence of solutions for
equation~\eqref{1} is used in Bégout and D\'iaz~\cite{MR3190983}.

\section{Main results}
\label{main}

Here, we state the main results of this paper.

\begin{thm}[\textbf{Existence}]
\label{thmexista1}
Let $\O$ an open subset of $\R^N$ be such that $|\O|<\infty$ and assume $0<m<1,$ $(a,b)\in\C^2$ and $F\in L^2(\O).$ If $\Re(b)<0$ then assume further that $\Im(b)\neq0$ or $-\frac1{\CP^2}<\Re(b),$ where $\CP$ is the Poincaré's constant in~\eqref{poincare} below. Then there exists at least a solution $u\in H^1_0(\O)$ of~\eqref{1}. In addition, Symmetry Property~$\ref{sym}$ below holds.
\end{thm}

\begin{sym}
\label{sym}
If furthermore, for any $\vR\in SO_N(\R),$ $\vR\O=\O$ and if $F$ is spherically symmetric then we may construct a solution which is additionally spherically symmetric. For $N=1,$ this means that if $F$ is an even $($respectively, an odd$)$ function then $u$ is also an even $($respectively, an odd$)$ function.
\end{sym}

\begin{thm}[\textbf{\textit{A priori} bound}]
\label{thmbound1}
Let $\O$ an open subset of $\R^N$ be such that $|\O|<\infty$ and assume $0<m<1,$ $(a,b)\in\C^2$ and $F\in L^2(\O).$ If $\Re(b)<0$ then assume further that $\Im(b)\neq0$ or $-\frac1{\CP^2}<\Re(b),$ where $\CP$ is the constant in~\eqref{poincare} below. Let $u\in H^1_0(\O)$ be any solution
to~\eqref{1}. Then we have the following estimate.
\begin{gather*}
\|u\|_{H^1_0(\O)}\le C,
\end{gather*}
where $C=C(\|F\|_{L^2(\O)},|\O|,|a|,|b|,N,m).$
\end{thm}

\begin{thm}[\textbf{Existence}]
\label{thmexista2}
Let $\O\subseteq\R^N$ be an open subset and assume $V\in L^\infty(\O;\R),$ $0<m<1,$ $(a,b,c)\in\C^3$ is such that $\Im(a)\le0,$ $\Im(b)<0$ and
$\Im(c)\le0.$ If $\Re(a)\le0$ then assume further that $\Im(a)<0.$ Then we have the following result.
\begin{enumerate}
 \item[]
  \begin{enumerate}[$1)$]
   \item
    For any $F\in L^2(\O),$ there exists at least a solution $u\in H^1_0(\O)\cap L^{m+1}(\O)$ to \eqref{2}.
   \item
    If we assume furthermore that $\O$ is bounded with a $C^1$ boundary then the conclusion~$1)$ still holds true with $u\in H^1(\O)$ and the
    boundary condition \eqref{neu} instead of $u\in H^1_0(\O).$
  \end{enumerate}
\end{enumerate}
If, in addition, $V$ is spherically symmetric then Symmetry Property~$\ref{sym}$ holds.
\end{thm}

\begin{rmk}
\label{rmkdefsol}
Here are some comments about boundary condition.
\begin{enumerate}
 \item[]
  \begin{enumerate}[1)]
   \item
	If $u\not\in C(\ovl\O)$ and $\O$ has not a $C^{0,1}$ boundary, the condition $u_{|\Gamma}=0$ does not make sense (in the sense of the trace) and,
	in this case, has to be understood as $u\in H^1_0(\O).$
   \item
	Assume that $\O$ is bounded and has a $C^{1,1}$ boundary. Let $u\in H^1(\O)$ be any solution to \eqref{2} with the boundary
	condition~\eqref{neu}. Then $u\in H^2(\O)$ and boundary condition $\frac{\partial u}{\partial\nu}_{|\Gamma}=0$ makes sense in the sense of the
	trace $\gamma\big(\nabla u.\nu\big)=0.$ If, in addition, $u\in C^1(\ovl\O)$ then obviously for any $x\in\Gamma,$ $\frac{\partial u}{\partial\nu}(x)=0.$
	Indeed, since $u\in H^1(\O),$ $\Delta u\in L^2(\O)$ and \eqref{2} makes sense almost everywhere in $\O,$ we have
	 $\gamma\left(\frac{\partial u}{\partial\nu}\right)\in H^{-\frac12}(\Gamma)$ and by Green's formula,
	\begin{multline}
		\Re\vint_\O\nabla u(x).\ovl{\nabla v(x)}\d x
			-\left\langle\gamma\left(\frac{\partial u}{\partial\nu}\right),\gamma(v)\right\rangle_{H^{-\frac12}(\Gamma),H^\frac12(\Gamma)}	\\
		+\Re\vint_\O f\big(u(x)\big)\ovl{v(x)}\d x=\Re\vint_\O F(x)\ovl{v(x)}\d x,
	\end{multline}
	for any $v\in H^1(\O),$ where $f(u)=a|u|^{-(1-m)}u+bu+cV^2u$ (see Lemma~4.1, Theorem~4.2 and Corollary~4.1, p.155, in Lions and
	Magenes~\cite{MR0146525} and (1,5,3,10) in Grisvard~\cite{MR3396210}, p.62). This implies that
	\begin{gather}
		\label{rmkdefsolequ}
			\left\langle\gamma\left(\frac{\partial u}{\partial\nu}\right),\gamma(v)\right\rangle_{H^{-\frac12}(\Gamma),H^\frac12(\Gamma)}=0,
	\end{gather}
	for any $v\in H^1(\O).$ Let $w\in H^\frac12(\Gamma).$ Let $v\in H^1(\O)$ be such that $\gamma(v)=w$ (Theorem~1.5.1.3, p.38, in
	Grisvard~\cite{MR3396210}). We then deduce from~\eqref{rmkdefsolequ} that,
	\begin{gather*}
		\forall w\in H^\frac12(\Gamma), \;
		\left\langle\gamma\left(\frac{\partial u}{\partial\nu}\right),w\right\rangle_{H^{-\frac12}(\Gamma),H^\frac12(\Gamma)}=0,
	\end{gather*}
	and so $\gamma\left(\frac{\partial u}{\partial\nu}\right)=0.$ But also $u\in L^2(\O)$ and $\Delta u\in L^2(\O).$ It follows that $u\in H^2(\O)$
	(Proposition~2.5.2.3, p.131, in Grisvard~\cite{MR3396210}). Hence the result.
  \end{enumerate}
\end{enumerate}
\end{rmk}

\begin{thm}[\textbf{\textit{A priori} bound}]
\label{thmbound2}
Let $\O\subseteq\R^N$ be an open subset, let $V\in L^\infty(\O;\R),$ let $0<m<1,$ let $(a,b,c)\in\C^3$ be such that $\Im(a)\le0,$ $\Im(b)<0$ and $\Im(c)\le0.$ If $\Re(a)\le0$ then assume further that $\Im(a)<0.$ Let $F\in L^2(\O)$ and let $u\in H^1(\O)$ be any solution to~\eqref{2} with boundary condition~\eqref{dir} or~\eqref{neu}\footnote{\label{fn1}for which we additionally assume that $\O$ has a $C^1$ boundary.}. Then we have the following estimate.
\begin{gather*}
\|u\|_{H^1(\O)}^2+\|u\|_{L^{m+1}(\O)}^{m+1}\le M\big(\|V\|_{L^\infty(\O)}^4+1\big)\|F\|_{L^2(\O)}^2,
\end{gather*}
where $M=M(|a|,|b|,|c|).$
\end{thm}

\begin{com}
\label{com}
In the context of the paper of Bégout and D\'iaz~\cite{MR2876246}, we can establish an existence result with the homogeneous Neumann boundary condition (instead of the homogeneous Dirichlet condition) and $F\in L^2(\O)$ \big(instead of $F\in L^\frac{m+1}{m}(\O)\big).$ In Bégout and
D\'iaz~\cite{MR2876246}, we introduced the set,
\begin{gather*}
\wt\A=\C\setminus\big\{z\in\C; \Re(z)=0 \text{ and } \Im(z)\le0\big\},
\end{gather*}
and assumed that $(\wt a,\wt b)\in\C^2$ satisfies,
\begin{gather}
\label{abtilde}
(\wt a,\wt b)\in\wt\A\times\wt\A \quad \text{ and } \quad
	\begin{cases}
		\Re(\wt a)\Re(\wt b)\ge0,													\medskip \\
		\text{ or }																\medskip \\
		\Re(\wt a)\Re(\wt b)<0 \; \text{ and } \; \Im(\wt b)>\dfrac{\Re(\wt b)}{\Re(\wt a)}\Im(\wt a),
	\end{cases}
\end{gather}
with possibly $\wt b=0,$ and we worked with
\begin{gather*}
-\vi\Delta u+\wt a|u|^{-(1-m)}u+\wt bu=\wt F.
\end{gather*}
Nevertheless, to maintain a closer notation to many applied works in the literature (see, e.g., the introduction of Bégout and D\'iaz~\cite{MR3193996}), we do not work any more with this equation but with,
\begin{gather*}
-\Delta u+a|u|^{-(1-m)}u+bu=F,
\end{gather*}
and $b\neq0.$ This means that we chose, $\wt a=\vi a,$ $\wt b=\vi b$ and $\wt F=\vi F.$ Then assumptions on $(a,b)$ are changed by the fact that for $\wt z=\vi z$,
\begin{gather}
\label{sca1}
\Re(z)=\Re(-\vi\wt z)=\Im(\wt z),	\\
\label{sca2}
\Im(z)=\Im(-\vi\wt z)=-\Re(\wt z).
\end{gather}
It follows that the set $\wt\A$ and \eqref{abtilde} become,
\begin{gather}
\label{A}
\A=\C\setminus\big\{z\in\C; \Re(z)\le0 \text{ and } \Im(z)=0\big\},	
\end{gather}
\begin{gather}
\label{ab}
(a,b)\in\A\times\A \quad \text{ and } \quad
	\begin{cases}
		\Im(a)\Im(b)\ge0,										\medskip \\
		\text{ or }												\medskip \\
		\Im(a)\Im(b)<0 \; \text{ and } \; \Re(b)>\dfrac{\Im(b)}{\Im(a)}\Re(a).
	\end{cases}
\end{gather}
Obviously,
\begin{gather*}
\Big((\wt a,\wt b)\in\wt\A\times\wt\A \text{ satisfies }\eqref{abtilde}\Big) \iff \Big((a,b)\in\A\times\A \text{ satisfies }\eqref{ab}\Big).
\end{gather*}
\end{com}

\noindent
Assumptions~\eqref{ab} are made to prove the existence and the localization property of solutions to equation~\eqref{1}. Now, we give some results about equation \eqref{1} when $(a,b)\in\A\times\A$ satisfies \eqref{ab}.

\begin{thm}[\textbf{Existence}]
\label{thmexista3}
Let $\O\subseteq\R^N$ be an open subset of $\R^N,$ let $0<m<1$ and let $(a,b)\in\A^2$ satisfies~\eqref{ab}.
\begin{enumerate}
 \item[]
  \begin{enumerate}[$1)$]
   \item
    For any $F\in L^2(\O),$ there exists at least a solution $u\in H^1_0(\O)\cap L^{m+1}(\O)$ to
	\begin{gather}
		\label{nls}
			-\Delta u+a|u|^{-(1-m)}u+bu=F, \text{ in } L^2(\O)+L^\frac{m+1}{m}(\O).
	\end{gather}
   \item
    If we assume furthermore that $\O$ is bounded with a $C^1$ boundary then the conclusion~$1)$ still holds true with $u\in H^1(\O)$ and the
    boundary condition \eqref{neu} instead of $u\in H^1_0(\O).$
  \end{enumerate}
\end{enumerate}
In addition, Symmetry Property~$\ref{sym}$ holds.
\end{thm}

\begin{thm}[\textbf{\textit{A priori} bound}]
\label{thmbound3}
Let $\O\subseteq\R^N$ be an open subset of $\R^N,$ let $0<m<1$ and let $(a,b)\in\A^2$ satisfies~\eqref{ab}. Let $F\in L^2(\O)$ and let
$u\in H^1(\O)\cap L^{m+1}(\O)$ be any solution to~\eqref{nls} with boundary condition~\eqref{dir} or~\eqref{neu}\textsuperscript{\textnormal{\ref{fn1}}}. Then we have the following estimate.
\begin{gather*}
\|u\|_{H^1(\O)}^2+\|u\|_{L^{m+1}(\O)}^{m+1}\le M\|F\|_{L^2(\O)}^2,
\end{gather*}
where $M=M(|a|,|b|).$
\end{thm}

\begin{thm}[\textbf{Uniqueness}]
\label{thmuni}
Let $\O\subseteq\R^N$ be an open subset, let $V\in L^\infty_\loc(\O;\R),$ let $0<m<1$ and let $(a,b,c)\in\C^3$ satisfies one of the three following conditions.
\begin{enumerate}
 \item[]
  \begin{enumerate}[$1)$]
   \item
    $a\neq0,$ $\Re(a)\ge0,$ $\Re(a\ovl b)\ge0$ and $\Re(a\ovl c)\ge0.$
   \item
    $b\neq0,$ $\Re(b)\ge0,$ $a=kb,$ for some $k\ge0$ and $\Re(b\ovl c)\ge0.$
   \item
    $c\neq0,$ $\Re(c)\ge0,$ $a=kc,$ for some $k>0$ and $\Re(b\ovl c)\ge0.$
  \end{enumerate}
\end{enumerate}
Let $F\in L^1_\loc(\O).$ If there exist two solutions  $u_1,u_2\in H^1(\O)\cap L^{m+1}(\O)$ of~\eqref{2} with the same boundary
condition~\eqref{dir} or~\eqref{neu}\textsuperscript{\textnormal{\ref{fn1}}} such that $Vu_1,Vu_2\in L^2(\O)$ then $u_1=u_2.$
\end{thm}

\begin{rmk}
\label{rmkthmexiuni}
Here are some comments about Theorems~\ref{thmexista1}, \ref{thmexista2}, \ref{thmexista3} and \ref{thmuni}.
\begin{enumerate}
 \item[]
  \begin{enumerate}[1)]
  \item
   Assume $F$ is spherically symmetric. Since we do not know, in general, if we have uniqueness of the solution, we are not able to show that any
   solution is radially symmetric.
   \item
    In Theorem~5.2 in Bégout and D\'iaz~\cite{MR2876246}, uniqueness for equation
	\begin{gather*}
		-\vi\Delta u+\wt a|u|^{-(1-m)}u+\wt bu=\wt F,
	\end{gather*}
    holds if $\wt a\neq0,$ $\Im(\wt a)\ge 0$ and $\Re(\wt a\ovl{\wt b})\ge0.$ By \eqref{sca1}--\eqref{sca2}, those assumptions are equivalent to~1) of
    Theorem~\ref{thmuni} above for equation \eqref{1} (of course, $c=0).$ It follows that Theorem~\ref{thmuni} above extends Theorem~5.2 of
    Bégout and D\'iaz~\cite{MR2876246}.
   \item
    In 2) of the above theorem, if we want to make an analogy with 1), assumption $a=kb,$ for some $k\ge0$ has to be replaced with $\Re(a\ovl b)\ge0$
    and $\Im(a\ovl b)=0.$ But,
	\begin{gather*}
		\Big(\Re(a\ovl b)\ge0 \text{ and } \Im(a\ovl b)=0\Big) \iff \Big(\exists k\ge0 \slash a=kb\Big).
	\end{gather*}
    In the same way,
	\begin{gather*}
		\Big(\Re(a\ovl c)> \text{ and } \Im(a\ovl c)=0\Big) \iff \Big(\exists k>0/a=kc\Big).
	\end{gather*}
   \item
    In the case of real solutions (with $F\equiv0$ and $(a,b,c)\in\R\times\R\times\{0\}),$ it is well-known that if $b<0$ then it may appear multiplicity of
    solutions (once $m\in(0,1)$ and $a>0).$ For more details, see Theorem~1 in D\'iaz and Hern\'andez~\cite{dh}.
  \end{enumerate}
\end{enumerate}
\end{rmk}

\begin{thm}[\textbf{Regularity}]
\label{thmreg}
Let $\O\subseteq\R^N$ be an open subset, let $V\in L^r_\loc(\O;\C),$ for any $1<r<\infty,$ let $0<m<1,$ let $a\in\C,$ let
$F\in L^1_\loc(\O),$ let $1<q<\infty$ and let $u\in L^q_\loc(\O)$ be any local solution to
\begin{gather}
\label{thmreg1}
-\Delta u+a|u|^{-(1-m)}u+Vu=F, \text{ in } \Dr^\p(\O).
\end{gather}
Let $q\le p<\infty$ and let $\alpha\in(0,m].$
\begin{enumerate}
 \item[]
  \begin{enumerate}[$1)$]
   \item
    If $F\in L^p_\loc(\O)$ then $u\in W^{2,p}_\loc(\O).$ If $(F,V)\in C_\loc^{0, \alpha}(\O)\times C_\loc^{0, \alpha}(\O)$ then $u\in C_\loc^{2,\alpha}(\O).$
   \item
    Assume further that $\O$ is bounded with a $C^{1,1}$ boundary, $F\in L^p(\O),$ $V\in L^r(\O;\C),$ for any $1<r<\infty,$ $u\in L^q(\O)$ and
    $\gamma(u)=0.$ Then $u\in W^{2,p}(\O)\cap W^{1,p}_0(\O).$ If $(F,V)\in C^{0,\alpha}(\ovl\O)\times C^{0,\alpha}(\ovl\O)$ then
    $u\in C^{2,\alpha}(\ovl\O)\cap C_0(\O).$
   \item
    Assume further that $\O$ is bounded with a $C^{1,1}$ boundary, $F\in L^p(\O),$ $V\in L^r(\O;\C),$ for any $1<r<\infty,$ $u\in L^q(\O)$ and
    $\gamma\left(\frac{\partial u}{\partial\nu}\right)=0.$ Then $u\in W^{2,p}(\O).$ If $(F,V)\in C^{0,\alpha}(\ovl\O)\times C^{0,\alpha}(\ovl\O)$ then
    $u\in C^{2,\alpha}(\ovl\O)$ and for any $x\in\Gamma,$ $\frac{\partial u}{\partial\nu}(x)=0.$
  \end{enumerate}
\end{enumerate}
\end{thm}

\begin{rmk}
\label{rmkthmreg}
Assume $\O$ is bounded and has a $C^{1,1}$ boundary. Let $V\in\bigcap\limits_{1<r<\infty}L^r(\O;\C),$ $0<m<1,$ $(a,b)\in\C^2,$ $1<q\le p<\infty,$
$F\in L^p(\O)$ and let $u\in L^q(\O)$ be any solution to~\eqref{thmreg1}. Let
$T:u\tends\left\{\gamma(u),\gamma\left(\frac{\partial u}{\partial\nu}\right)\right\}$ be the trace function defined on
$\Dr(\ovl\Omega).$ By density of $\Dr(\ovl\Omega)$ in
$D_q(\Delta)\stackrel{\mathrm{def}}{=}\big\{u\in L^q(\Omega);\Delta u\in L^q(\Omega)\big\},$ $T$ has a linear and continuous extension from $D_q(\Delta)$ into $W^{-\frac1q,q}(\Gamma)\times W^{-1-\frac1q,q}(\Gamma)$ (Hörmander~\cite{MR0106333},
Theorem~2 p.503; Lions and Magenes~\cite{MR0146525}, Lemma~2.2 and Theorem~2.1 p.147; Lions and Magenes~\cite{MR0146526}, Propositions~9.1,~9.2 and Theorem~9.1 p.82; Grisvard~\cite{MR3396210}, p.54). Since $u\in L^q(\Omega),$ it follows from equation~\eqref{thmreg1} and Hölder's inequality that $u\in D_q(\Delta),$ so that ``$\gamma(u)=0$'' and ``$\gamma\left(\frac{\partial u}{\partial\nu}\right)=0$'' make sense.
\end{rmk}

\noindent
The main difficulty to apply Theorem~\ref{thmreg} is to show that such a solution of~\eqref{thmreg1} verifies some boundary condition. In the following result, we give a sufficient condition.

\begin{prop}[\textbf{Regularity}]
\label{propreg}
Let $\O$ be a bounded open subset of $\R^N$ with a $C^{1,1}$ boundary, let $V\in L^N(\O;\C)$ $(V\in L^{2+\eps}(\O;\C),$ for some $\eps>0,$ if $N=2$ and $V\in L^2(\O;\C)$ if $N=1),$ let $0<m<1,$ let $a\in\C$ and let $F\in L^2(\O).$
\begin{enumerate}
 \item[]
  \begin{enumerate}[$1)$]
   \item
    Let $u\in H^1_0(\O)$ be any solution to \eqref{thmreg1}. Then $u\in H^2(\O)$ and $\gamma(u)=0.$
   \item
    Let $u\in H^1(\O)$ be any solution to \eqref{thmreg1} and \eqref{neu}. Then $u\in H^2(\O)$ and $\gamma\left(\frac{\partial u}{\partial\nu}\right)=0.$
  \end{enumerate}
\end{enumerate}
\end{prop}

\begin{rmk}
\label{rmkthmreg1}
Any solution given by Theorems~\ref{thmexista1}, \ref{thmexista2} or \ref{thmexista3} belongs to $H^2_\loc(\O)$ (Theorem~\ref{thmreg}).
\end{rmk}

\section{Notations}
\label{not}

We indicate here some of the notations used throughout this paper which have not been defined yet in the introduction (Section~\ref{intro}). We write $\vi^2=-1.$ We denote by $\ovl z$ the conjugate of the complex number $z,$ $\Re(z)$ its real part and $\Im(z)$ its imaginary part. For $1\le p\le\infty,$ $p^\p$ is the conjugate of $p$ defined by $\frac1p+\frac1{p^\p}=1.$ The symbol $\O$ always indicates a nonempty open subset of $\R^N$ (bounded or not); its closure is denoted by $\ovl\O$ and its boundary by $\Gamma.$ For $A\in\{\O;\ovl\O\},$ the space $C(A)=C^0(A)$ is the set of continuous functions from $A$ to $\C$ and $C^k(A)$ $(k\in\N)$ is the space of functions lying in $C(A)$ and having all derivatives of order lesser or equal than $k$ belonging to $C(A).$ For $0<\alpha\le1$ and $k\in\N_0\stackrel{\mathrm{def}}{=}\N\cup\{0\},$ $C^{k,\alpha}_\loc(\Omega)
=\left\{u\in C^k(\Omega);\forall\omega\Subset\Omega,\;\vsum_{|\beta|=k}H^\alpha_\omega(D^\beta u)<+\infty\right\},$ where
$H^\alpha_\omega(u)=\vsup_{\left\{\substack{(x,y)\in\omega^2 \\ x\not=y}\right.}\frac{|u(x)-u(y)|}{|x-y|^\alpha}.$ The notation $\omega\Subset\O$ means that $\omega$ is a bounded open subset of $\R^N$ and $\ovl\omega\subset\O.$ In the same way, $C^{k,\alpha}(\ovl\Omega)=\left\{u\in C^k(\ovl\Omega);\vsum_{|\beta|=k}H^\alpha_\O(D^\beta u)<+\infty\right\}.$ The space $C_0(\O)$ consists of functions belonging to $C(\ovl\O)$ and vanishing at the boundary $\Gamma,$ $\Dr(\Omega)$ is the space of $C^\infty$ functions with compact support and $\Dr(\ovl\O)$ is the restriction to $\ovl\O$ of functions lying in $\Dr(\R^N).$ The trace function defined on $\Dr(\ovl\O)$ is denoted by $\gamma.$ For $1\le p\le\infty$ and $m\in\N,$ the usual Lebesgue and Sobolev spaces are respectively denoted by $L^p(\Omega)$ and $W^{m,p}(\Omega),$ $W^{m,p}_0(\Omega)$ is the closure of $\Dr(\Omega)$ under the $W^{m,p}$-norm, $H^m(\Omega)=W^{m,2}(\Omega)$ and $H^m_0(\Omega)=W^{m,2}_0(\Omega).$ For a Banach space $E,$ its topological dual is denoted by $E^\star$ and $\langle\: . \; , \: . \:\rangle_{E^\star,E}\in\R$ is the $E^\star-E$ duality product. In particular, for any $T\in L^{p^\p}(\O)$ and $\vphi\in L^p(\O)$ with $1\le p<\infty,$ $\langle T,\vphi\rangle_{L^{p^\p}(\O),L^p(\O)}=\Re\vint_\O T(x)\ovl{\vphi(x)}\d x.$ We write, $W^{-m,p^\p}(\Omega)=\left(W^{m,p}_0(\Omega)\right)^\star$ $(p<\infty)$ and $H^{-m}(\Omega)=\left(H^m_0(\Omega)\right)^\star.$ Unless if specified, any function belonging in a functional space $\big(W^{m,p}(\Omega),$ $C^k(\Omega),$ etc\big) is supposed to be a complex-valued function $\big(W^{m,p}(\Omega;\C),$ $C^k(\Omega;\C),$ etc\big). We denote by $SO_N(\R)$ the special orthogonal group of $\R^N.$ Finally, we denote by $C$ auxiliary positive constants, and sometimes, for positive parameters $a_1,\ldots,a_n,$ write $C(a_1,\ldots,a_n)$ to indicate that the constant $C$ continuously depends only on $a_1,\ldots,a_n$ (this convention also holds for constants which are not denoted by ``$C$'').

\section{A priori estimates}
\label{apriori}

The proofs of the existence theorems relies on \textit{a priori} bounds, in order to truncate the nonlinearity and pass to the limit. These bounds are formally obtained by multiplying the equation by $\ovl u$ and $\vi\ovl u,$ integrate by parts and by making some linear combinations with the obtained results. Now, we recall the well-known Poincaré's inequality. If $|\O|<\infty$ then,
\begin{gather}
 \label{poincare}
  \forall u\in H^1_0(\O), \; \|u\|_{L^2(\O)}\le\CP\|\nabla u\|_{L^2(\O)}.
\end{gather}
where $\CP=\CP(|\O|,N).$ We will frequently use Hölder's inequality in the following form. If $|\O|<\infty$ and $0\le m\le1$ then
$L^2(\O)\inj L^{m+1}(\O)$ and
\begin{gather}
 \label{Lm}
  \forall u\in L^2(\O), \; \|u\|_{L^{m+1}(\O)}^{m+1}\le|\O|^\frac{1-m}{2}\|u\|_{L^2(\O)}^{m+1}.
\end{gather}

\noindent
Finally, we recall the well-known Young's inequality. For any real $x\ge0,$ $y\ge0$ and $\mu>0,$ one has
\begin{gather}
\label{young}
xy\le\frac{\mu^2}2x^2+\frac1{2\mu^2}y^2.
\end{gather}

\begin{lem}
\label{lemeap1}
Let $\O$ an open subset of $\R^N$ be such that $|\O|<\infty,$ let $\o$ an open subset of $\R^N$ be such that $\o\subseteq\O,$ let $0\le m\le 1,$ let
$(a,b)\in\C^2,$ let $\alpha,\beta\ge0$ and let $F\in L^2(\O).$ Let $u\in H^1_0(\O)$ satisfies
\begin{multline}
\label{lemeap11}
\Big|\|\nabla u\|_{L^2(\O)}^2+\Re(a)\left(\|u\|_{L^{m+1}(\o)}^{m+1}+\alpha\|u\|_{L^1(\o^\c)}\right)\\
+\Re(b)\left(\|u\|_{L^2(\o)}^2+\beta\|u\|_{L^1(\o^\c)}\right)\Big|\le\int_\O|Fu|\d x,
\end{multline}
\begin{gather}
\label{lemeap12}
\left|\Im(a)\left(\|u\|_{L^{m+1}(\o)}^{m+1}+\alpha\|u\|_{L^1(\o^\c)}\right)+\Im(b)\left(\|u\|_{L^2(\o)}^2+\beta\|u\|_{L^1(\o^\c)}\right)\right|\le\int_\O|Fu|\d x.
\end{gather}
Here, $\o^\c=\O\setminus\o.$ Assume that one of the three following assertions holds.
\begin{enumerate}
	\item[]
		\begin{enumerate}[$1)$]
			\item
				$\Re(b)\ge0.$ If $\Re(a)<0$ and $|\o|<|\O|$ then assume further that $\alpha\|u\|_{L^1(\o^\c)}\le\|u\|_{L^{m+1}(\o^\c)}^{m+1}.$
			\item
				$\Re(b)<0$ and $\Im(b)\neq0.$ If $|\o|<|\O|$ then assume further that $\alpha\|u\|_{L^1(\o^\c)}\le\|u\|_{L^{m+1}(\o^\c)}^{m+1},$
				$F\in L^\infty(\O)$ and $-\alpha|\Im(a)|+\frac\beta2|\Im(b)|>\|F\|_{L^\infty(\O)}.$
			\item
				$-\CP^{-2}<\Re(b)<0,$ where $\CP$ is the constant in~\eqref{poincare},
				$\alpha\|u\|_{L^1(\o^\c)}\le\|u\|_{L^{m+1}(\o^\c)}^{m+1}$ and $\beta\|u\|_{L^1(\o^\c)}\le\|u\|_{L^2(\o^\c)}^2.$
		\end{enumerate}
\end{enumerate}
Then we have the following estimate.
\begin{gather}
\label{eap1}
\|u\|_{H^1_0(\O)}\le C,
\end{gather}
where $C=C(\|F\|_{L^2(\O)},|\O|,|a|,|b|,N,m).$
\end{lem}

\begin{rmk}
\label{rmklemeap1}
Obviously, if $|\o|=|\O|$ then $\alpha\|u\|_{L^1(\o^\c)}\le\|u\|_{L^{m+1}(\o^\c)}^{m+1}$ and $\beta\|u\|_{L^1(\o^\c)}\le\|u\|_{L^2(\o^\c)}^2.$
\end{rmk}

\begin{vproof}{of Lemma~\ref{lemeap1}.}
By Poincaré's inequality~\eqref{poincare}, it is sufficient to establish
\begin{gather}
\label{prooflemeap1}
\|\nabla u\|_{L^2(\O)}\le C(\|F\|_{L^2(\O)},|\O|,|a|,|b|,N,m).
\end{gather}
Moreover, it follows from \eqref{young} and \eqref{poincare},
\begin{gather}
\label{prooflemeap2}
\int_\O|Fu|\d x\le\frac{\CP^2}2\|F\|_{L^2(\O)}^2+\frac12\|\nabla u\|_{L^2(\O)}^2.
\end{gather}
Finally, it follows from \eqref{Lm} and \eqref{poincare} that if $\alpha\|u\|_{L^1(\o^\c)}\le\|u\|_{L^{m+1}(\o^\c)}^{m+1}$ then one has,
\begin{gather}
\label{prooflemeap3}
\|u\|_{L^{m+1}(\o)}^{m+1}+\alpha\|u\|_{L^1(\o^\c)}\le\|u\|_{L^{m+1}(\O)}^{m+1}\le\CP^{m+1}|\O|^\frac{1-m}{2}\|\nabla u\|_{L^2(\O)}^{m+1}.
\end{gather}
We divide the proof in 3 steps. \\
\textbf{Step 1.} Proof of~\eqref{prooflemeap1} with Assumption~1). \\
Assume hypothesis~1) holds true. If $\Re(a)\ge0$ then \eqref{prooflemeap1} follows from \eqref{lemeap11} and \eqref{prooflemeap2}, while if $\Re(a)<0$ we then deduce from  \eqref{lemeap11}, \eqref{prooflemeap2} and \eqref{prooflemeap3} that,
\begin{gather*}
\left(\|\nabla u\|_{L^2(\O)}^{1-m}-2|\Re(a)|\CP^{m+1}|\O|^\frac{1-m}{2}\right)\|\nabla u\|_{L^2(\O)}^{m+1}\le\CP^2\|F\|_{L^2(\O)}^2.
\end{gather*}
Hence \eqref{prooflemeap1}. \\
\textbf{Step 2.} Proof of~\eqref{prooflemeap1} with Assumption~2). \\
As for Step~1, it follows from \eqref{lemeap12}, \eqref{Lm}, \eqref{young} and Hölder's inequality that
\begin{multline*}
|\Im(b)|\left(\|u\|_{L^2(\o)}^2+\beta\|u\|_{L^1(\o^\c)}\right)	\le |\Im(a)||\O|^\frac{1-m}{2}\|u\|_{L^2(\o)}^{m+1}+\alpha|\Im(a)|\|u\|_{L^1(\o^\c)}	 \\
+\frac1{2|\Im(b)|}\|F\|_{L^2(\o)}^2+\frac{|\Im(b)|}2\|u\|_{L^2(\o)}^2+\|F\|_{L^\infty(\o^\c)}\|u\|_{L^1(\o^\c)}.
\end{multline*}
Recalling that when $|\o|<|\O|,$ $-\alpha|\Im(a)|+\frac\beta2|\Im(b)|>\|F\|_{L^\infty(\O)},$ the above estimate yields
\begin{gather}
\label{prooflemeap4}
\left(|\Im(b)|\|u\|_{L^2(\o)}^{1-m}-2|\Im(a)||\O|^\frac{1-m}{2}\right)\|u\|_{L^2(\o)}^{m+1}+\beta|\Im(b)|\|u\|_{L^1(\o^\c)}\le\frac1{|\Im(b)|}\|F\|_{L^2(\o)}^2.
\end{gather}
If $|\Im(b)|\|u\|_{L^2(\o)}^{1-m}-2|\Im(a)||\O|^\frac{1-m}{2}\le1$ then 
\begin{gather}
\label{prooflemeap5}
\|u\|_{L^2(\o)}\le C(\|F\|_{L^2(\O)},|\O|,|a|,|b|,m)\stackrel{\mathrm{not.}}{=}C_0,
\end{gather}
and it follows from \eqref{lemeap12}, \eqref{Lm}, \eqref{prooflemeap5} and Hölder's inequality that,
\begin{multline*}
\big(\beta|\Im(b)|-\alpha|\Im(a)|\big)\|u\|_{L^1(\o^\c)}\le C(C_0)+\|F\|_{L^\infty(\o^\c)}\|u\|_{L^1(\o^\c)}		\\
\le C(C_0)+\left(\frac\beta2|\Im(b)|-\alpha|\Im(a)|\right)\|u\|_{L^1(\o^\c)},
\end{multline*}
so that,
\begin{gather}
\label{prooflemeap6}
\beta\|u\|_{L^1(\o^\c)}\le C(\|F\|_{L^2(\O)},|\O|,|a|,|b|,m)\stackrel{\mathrm{not.}}{=}C_1.
\end{gather}
But if $|\Im(b)|\|u\|_{L^2(\o)}^{1-m}-2|\Im(a)||\O|^\frac{1-m}{2}>1$ then \eqref{prooflemeap5} and \eqref{prooflemeap6} come from \eqref{prooflemeap4}. \\
Finally, by~\eqref{lemeap11}, \eqref{prooflemeap2}, \eqref{prooflemeap3}, \eqref{prooflemeap5} and \eqref{prooflemeap6}, one obtains
\begin{gather*}
\|\nabla u\|_{L^2(\O)}^2
\le|\Re(a)|\CP^{m+1}|\O|^\frac{1-m}{2}\|\nabla u\|_{L^2(\O)}^{m+1}+C(C_0,C_1)+\frac{\CP^2}2\|F\|_{L^2(\O)}^2+\frac12\|\nabla u\|_{L^2(\O)}^2.
\end{gather*}
It follows that $\big(\|\nabla u\|_{L^2(\O)}^{1-m}-C\big)\|\nabla u\|_{L^2(\O)}^{m+1}\le C+\CP^2\|F\|_{L^2(\O)}^2,$ from which we easily deduce \eqref{prooflemeap1}. \\
\textbf{Step 3.} Proof of~\eqref{prooflemeap1} with Assumption~3). \\
Let $\mu>0.$ By Assumption~3), \eqref{poincare}, \eqref{young} and \eqref{prooflemeap3}
\begin{gather*}
\|\nabla u\|_{L^2(\O)}^2
\le C\|\nabla u\|_{L^2(\O)}^{m+1}+\left(|\Re(b)|\CP^2+\frac{\CP^2}{2\mu^2}\right)\|\nabla u\|_{L^2(\O)}^2+\frac{\mu^2}2\|F\|_{L^2(\O)}^2,
\end{gather*}
where $C=C(|\O|,|a|,N,m).$ We then deduce,
\begin{gather*}
\left(\left(1-|\Re(b)|\CP^2-\frac{\CP^2}{2\mu^2}\right)\|\nabla u\|_{L^2(\O)}^{1-m}-C\right)\|\nabla u\|_{L^2(\O)}^{m+1}\le\frac{\mu^2}2\|F\|_{L^2(\O)}^2.
\end{gather*}
Since $|\Re(b)|<\CP^{-2},$ there exists $\mu_0>0$ such that $C_2\stackrel{\mathrm{def.}}{=}1-|\Re(b)|\CP^2-\frac{\CP^2}{2\mu_0^2}>0.$
For such a $\mu_0,$ $\left(C_2\|\nabla u\|_{L^2(\O)}^{1-m}-C\right)\|\nabla u\|_{L^2(\O)}^{m+1}\le\frac{\mu_0^2}2\|F\|_{L^2(\O)}^2,$ from which \eqref{prooflemeap1} follows.
\medskip
\end{vproof}

\begin{cor}
\label{corbound1}
Let $(\O_n)_{n\in\N}$ a sequence of open subsets of $\R^N$ be such that $\vsup_{n\in\N}|\O_n|<\infty,$ let $0<m<1,$ let $(a,b)\in\C^2$ and let $(F_n)_{n\in\N}\subset L^\infty(\O_n)$ be such that $\vsup_{n\in\N}\|F_n\|_{L^2(\O_n)}<\infty.$ If $\Re(b)<0$ then assume further that $\Im(b)\neq0$ or $-\frac1{\CP^2}<\Re(b),$ where $\CP$ is the constant in~\eqref{poincare}. Let $(u^n_\ell)_{(n,\ell)\in\N^2}\subset H^1_0(\O_n)$ be a sequence satisfying
\begin{gather}
\label{corbound11}
\forall n\in\N, \; \forall\ell\in\N, \; -\Delta u_\ell^n+f_\ell\big(u_\ell^n\big)=F_n, \text{ in } L^2(\O_n),
\end{gather}
where for any $\ell\in\N,$
\begin{gather}
\label{corbound12}
\forall u\in L^2(\O_n), \; f_\ell(u)=
	\begin{cases}
		a|u|^{-(1-m)}u+bu,					&	\mbox{ if } |u|\le\ell, \medskip \\
		a\ell^m\dfrac{u}{|u|}+b\ell\dfrac{u}{|u|},	&	\mbox{ if } |u| > \ell.
	\end{cases}
\end{gather}
Then there exists a diagonal extraction $\left(u_{\vphi(n)}^n\right)_{n\in\N}$ of $(u_\ell^n)_{(n,\ell)\in\N^2}$ such that the following estimate holds.
\begin{gather*}
\forall n\in\N, \; \big\|u_{\vphi(n)}^n\big\|_{H^1_0(\O_n)}\le C,
\end{gather*}
where $C=C\left(\vsup_{n\in\N}\|F_n\|_{L^2(\O_n)},\vsup_{n\in\N}|\O_n|,|a|,|b|,N,m\right).$
\end{cor}

\begin{proof*}
Choosing $u_\ell^n$ and $\vi u_\ell^n$ as test functions, we get
\begin{multline*}
\|\nabla u_\ell^n\|_{L^2(\O_n)}^2+\Re(a)\left(\|u_\ell^n\|_{L^{m+1}(\{|u_\ell^n|\le\ell\})}^{m+1}+\ell^m\|u_\ell^n\|_{L^1(\{|u_\ell^n|>\ell\})}\right)		\\
	+\Re(b)\left(\|u_\ell^n\|_{L^2(\{|u_\ell^n|\le\ell\})}^2+\ell\|u_\ell^n\|_{L^1(\{|u_\ell^n|>\ell\})}\right)=\Re\int_{\O_n}F_n\ovl{u_\ell^n}\d x,
\end{multline*}
\begin{multline*}
\Im(a)\left(\|u_\ell^n\|_{L^{m+1}(\{|u_\ell^n|\le\ell\})}^{m+1}+\ell^m\|u_\ell^n\|_{L^1(\{|u_\ell^n|>\ell\})}\right)							\\
	+\Im(b)\left(\|u_\ell^n\|_{L^2(\{|u_\ell^n|\le\ell\})}^2+\ell\|u_\ell^n\|_{L^1(\{|u_\ell^n|>\ell\})}\right)=\Im\int_{\O_n}F_n\ovl{u_\ell^n}\d x,
\end{multline*}
for any $(n,\ell)\in\N^2.$ We first note that,
\begin{gather}
\label{proofcorbound11}
\forall(n,\ell)\in\N^2,
	\begin{cases}
		\ell^m\|u_\ell^n\|_{L^1(\{|u_\ell^n|>\ell\})}\le\|u_\ell^n\|_{L^{m+1}(\{|u_\ell^n|>\ell\})}^{m+1},	\medskip \\
		\ell\|u_\ell^n\|_{L^1(\{|u_\ell^n|>\ell\})}\le\|u_\ell^n\|_{L^2(\{|u_\ell^n|>\ell\})}^2,
	\end{cases}
\end{gather}
For each $n\in\N,$ we choose $\vphi(n)\in\N$ large enough to have $\vphi(n)^{1-m}>2\frac{\|F_n\|_{L^\infty(\O_n)}+|\Im(a)|}{|\Im(b)|},$ when
$\Im(b)\neq0$ and $\vphi(n)=n,$ when $\Im(b)=0.$ Thus for any $n\in\N,$ as soon as $\Im(b)\neq0,$ one has
\begin{gather}
\label{proofcorbound12}
\|F_n\|_{L^\infty(\O_n)}<-\vphi(n)^m|\Im(a)|+\dfrac{\vphi(n)}2|\Im(b)|.
\end{gather}
With help of \eqref{proofcorbound11} and \eqref{proofcorbound12}, we may apply Lemma~\ref{lemeap1} to $u_{\vphi(n)}^n,$ for each $n\in\N,$ with $\o=\left\{x\in\O_n;\left|u_{\vphi(n)}^n(x)\right|\le\vphi(n)\right\},$ $\alpha=\vphi(n)^m$ and $\beta=\vphi(n).$
\medskip
\end{proof*}

\begin{lem}
\label{lemeap2}
Let $\O\subseteq\R^N$ be an open subset, let $\o$ an open subset of $\R^N$ be such that $\o\subseteq\O,$ let $m\ge0$ and let $(a,b,c)\in\C^3$  be such that $\Im(b)\neq0.$ If $\Re(a)\le0$ then assume further that $\Im(a)\neq0.$ Let $\alpha,\beta,R\ge0,$ let $F\in L^2(\O)$ and let
\begin{gather*}
A=
	\begin{cases}
		\max\left\{1,\frac{1+|b|+R^2|c|}{|\Im(b)|},\frac{|\Re(a)|}{|\Im(a)|}\right\},		&	\text{if } \Re(a)\le0,	\medskip \\
		\max\left\{1,\frac{1+|b|+R^2|c|}{|\Im(b)|}\right\},						&	\text{if } \Re(a)>0.
	\end{cases}
\end{gather*}
If $|\o|<|\O|$ then assume further that $F\in L^\infty(\O)$ and $\beta\ge2A\|F\|_{L^\infty(\O)}+1.$ Let $u\in H^1(\O)$ satisfies
\begin{multline}
\label{lemeap2-1}
\|\nabla u\|_{L^2(\O)}^2+\Re(a)\left(\|u\|_{L^{m+1}(\o)}^{m+1}+\alpha\|u\|_{L^1(\o^\c)}\right) \\
-(|b|+R^2|c|)\left(\|u\|_{L^2(\o)}^2+\beta\|u\|_{L^1(\o^\c)}\right)\le\int_\O|Fu|\d x,
\end{multline}
\begin{gather}
\label{lemeap2-2}
|\Im(a)|\left(\|u\|_{L^{m+1}(\o)}^{m+1}+\alpha\|u\|_{L^1(\o^\c)}\right)+|\Im(b)|\left(\|u\|_{L^2(\o)}^2+\beta\|u\|_{L^1(\o^\c)}\right)\le\int_\O|Fu|\d x.
\end{gather}
Then there exists a positive constant $M=M(|a|,|b|,|c|)$ such that,
\begin{gather}
\label{eap2}
\|\nabla u\|_{L^2(\O)}^2+\|u\|_{L^2(\o)}^2+\|u\|_{L^{m+1}(\o)}^{m+1}+\|u\|_{L^1(\o^\c)}\le M(R^4+1)\|F\|_{L^2(\O)}^2.
\end{gather}
\end{lem}

\begin{proof*}
Let $A$ be as in the lemma. We multiply~\eqref{lemeap2-2} by $A$ and sum the result to~\eqref{lemeap2-1}. This yields,
\begin{gather*}
\|\nabla u\|_{L^2(\O)}^2+A_0\left(\|u\|_{L^{m+1}(\o)}^{m+1}+\alpha\|u\|_{L^1(\o^\c)}\right)
+\|u\|_{L^2(\o)}^2+\beta\|u\|_{L^1(\o^\c)}\le2A\int_\O|Fu|\d x,
\end{gather*}
where $A_0=A|\Im(a)|+\Re(a).$ Applying Hölder's inequality and \eqref{young}, we get
\begin{multline*}
\|\nabla u\|_{L^2(\O)}^2+\|u\|_{L^2(\o)}^2+A_0\|u\|_{L^{m+1}(\o)}^{m+1}+\beta\|u\|_{L^1(\o^\c)}	\\
\le2A\|F\|_{L^\infty(\O)}\|u\|_{L^1(\o^\c)}+2A^2\|F\|_{L^2(\O)}^2+\frac12\|u\|_{L^2(\o)}^2,
\end{multline*}
from which we deduce the result if $|\o|=|\O|.$ Now, suppose $|\o|<|\O|.$ The above estimate leads to,
\begin{gather*}
\|\nabla u\|_{L^2(\O)}^2+\|u\|_{L^2(\o)}^2+A_0\|u\|_{L^{m+1}(\o)}^{m+1}+\left(\beta-2A\|F\|_{L^\infty(\O)}\right)\|u\|_{L^1(\o^\c)}\le4A^2\|F\|_{L^2(\O)}^2,
\end{gather*}
from which we prove the lemma since $\beta-2A\|F\|_{L^\infty(\O)}\ge1.$
\medskip
\end{proof*}

\begin{lem}
\label{lemAB}
Let $(a,b)\in\A^2$ satisfies \eqref{ab}. Then there exists $\delta_\star=\delta_\star(|a|,|b|)\in(0,1],$ $L=L(|a|,|b|)$ and $M=M(|a|,|b|)$ satisfying the following property. If $\delta\in[0,\delta_\star]$ and $C_0,$ $C_1,$ $C_2,$ $C_3,$ $C_4$ are six nonnegative real numbers satisfying
\begin{gather}
\label{Re}
\big|C_1+\delta C_2+\Re(a)C_3+\big(\Re(b)-\delta\big)C_4\big|\le C_0,	\\
\label{Im}
\big|\Im(a)C_3+\Im(b)C_4\big|\le C_0,
\end{gather}
then
\begin{gather}
\label{lemestimAB}
0\le C_1+LC_3+LC_4\le MC_0.
\end{gather}
\end{lem}

\begin{proof*}
We split the proof in 4 cases. Let $\gamma>0$ be small enough to be chosen later. Note that when $\Im(a)\Im(b)\ge0$ then estimate~\eqref{Im} can be rewritten as
\begin{gather}
\label{Imm}
|\Im(a)|C_3+|\Im(b)|C_4\le C_0.
\end{gather}
\textbf{Case 1.} $\Re(a)\ge0,$ $\Re(b)\ge0$ and $\Im(a)\Im(b)\ge0.$ We add~\eqref{Imm} with \eqref{Re} and obtain,
\begin{gather*}
C_1+\big(\Re(a)+|\Im(a)|\big)C_3+\big(\Re(b)-\delta_\star+|\Im(b)|\big)C_4\le2C_0.
\end{gather*}
\textbf{Case 2.} $\Big(\Re(a)\ge0,$ $\Re(b)<0$ and $\Im(a)\Im(b)\ge0\Big)$ or $\Big(\Im(a)\Im(b)<0\Big).$ Then,
\begin{gather*}
C_1+\frac{\Re(a)\Im(b)-\Re(b)\Im(a)+\gamma\Im(a)}{\Im(b)}C_3+(\gamma-\delta_\star)C_4\le\frac{|\Re(b)|+|\Im(b)|+\gamma}{|\Im(b)|}C_0.
\end{gather*}
where we computed $\eqref{Re}-\frac{\Re(b)-\gamma}{\Im(b)}\eqref{Im}.$ \\
\textbf{Case 3.} $\Re(a)<0,$ $\Re(b)\ge0$ and $\Im(a)\Im(b)\ge0.$ By computing $\eqref{Re}-\frac{\Re(a)-\gamma}{\Im(a)}\eqref{Im},$ we get,
\begin{gather*}
C_1+\gamma C_3+\left(\frac{\Re(b)\Im(a)-\Re(a)\Im(b)+\gamma\Im(b)}{\Im(a)}-\delta_\star\right)C_4\le\frac{|\Re(a)|+|\Im(a)|+\gamma}{|\Im(a)|}C_0.
\end{gather*}
\textbf{Case 4.} $\Re(a)<0,$ $\Re(b)<0$ and $\Im(a)\Im(b)\ge0.$ Note that since $(a,b)\in\A^2$ then necessarily $\Im(a)\Im(b)\neq0.$ Thus, we can compute
$\eqref{Re}+\max\left\{\frac{|\Re(a)|+\gamma}{|\Im(a)|},\frac{|\Re(b)|+\gamma}{|\Im(b)|}\right\}\eqref{Imm}$ and obtain,
\begin{gather*}
C_1+\gamma C_3+(\gamma-\delta_\star)C_4\le\left(\frac{|\Re(a)|+|\Im(a)|+\gamma}{|\Im(a)|}+\frac{|\Re(b)|+|\Im(b)|+\gamma}{|\Im(b)|}\right)C_0.
\end{gather*}
In both cases, we may choose $\gamma>0$ small enough to have
\begin{gather*}
	\begin{cases}
		\dfrac{\Re(a)\Im(b)-\Re(b)\Im(a)+\gamma\Im(a)}{\Im(b)}>0,	&	\text{in Case~2}, \medskip \\
		\dfrac{\Re(b)\Im(a)-\Re(a)\Im(b)+\gamma\Im(b)}{\Im(a)}>0,	&	\text{in Case~3}.
	\end{cases}
\end{gather*}
Then we choose $0<\delta_\star<\min\big\{1,\gamma,|\Im(b)|+|\Re(b)|\big\}$ such that
\begin{gather*}
\delta_\star<\frac{\Re(b)\Im(a)-\Re(a)\Im(b)+\gamma\Im(b)}{\Im(a)}, \text{ in Case~3}.
\end{gather*}
This ends the proof.
\medskip
\end{proof*}

\begin{cor}
\label{corbound2}
Let $\O\subseteq\R^N$ be an open subset, let $V\in L^\infty(\O;\R),$ let $0<m<1$ and let $(a,b,c)\in\C^3$ be such that $\Im(a)\le0,$ $\Im(b)<0$ and
$\Im(c)\le0.$ If $\Re(a)\le0$ then assume further that $\Im(a)<0.$ Let $\delta\ge0.$ Let $(F_n)_{n\in\N}\subset L^\infty(\O)\cap L^2(\O)$ be bounded in $L^2(\O)$ and let $(u^n_\ell)_{(n,\ell)\in\N^2}\subset H^1(\O)\cap L^{m+1}(\O)$ be a sequence satisfying
\begin{gather}
\label{corbound2-1}
\forall n\in\N, \; \forall\ell\in\N, \; -\Delta u_\ell^n+\delta u_\ell^n+f_\ell\big(u_\ell^n\big)=F_n, \text{ in } L^2(\O),
\end{gather}
with boundary condition~\eqref{dir} or \eqref{neu}, where for any $\ell\in\N,$
\begin{gather}
\label{corbound2-2}
\forall u\in L^2(\O), \; f_\ell(u)=
	\begin{cases}
		a|u|^{-(1-m)}u+(b-\delta)u+cV^2u,								&	\mbox{ if } |u|\le\ell, \medskip \\
		a\ell^m\dfrac{u}{|u|}+(b-\delta)\ell\dfrac{u}{|u|}+cV^2\ell\dfrac{u}{|u|},		&	\mbox{ if } |u| > \ell.
	\end{cases}
\end{gather}
For \eqref{neu}, $\O$ is assumed to have a $C^1$ boundary. Then there exist $M=M\big(\|V\|_{L^\infty(\O)},|a|,|b|,|c|\big)$ and a diagonal extraction
$\left(u_{\vphi(n)}^n\right)_{n\in\N}$ of $(u_\ell^n)_{(n,\ell)\in\N^2}$ for which,
\begin{multline*}
\big\|\nabla u_{\vphi(n)}^n\big\|_{L^2(\O)}^2+\big\|u_{\vphi(n)}^n\big\|_{L^2\left(\left\{\left|u_{\vphi(n)}^n\right|\le\vphi(n)\right\}\right)}^2
+\big\|u_{\vphi(n)}^n\big\|_{L^{m+1}\left(\left\{\left|u_{\vphi(n)}^n\right|\le\vphi(n)\right\}\right)}^{m+1}								\\
+\big\|u_{\vphi(n)}^n\big\|_{L^1\left(\left\{\left|u_{\vphi(n)}^n\right|>\vphi(n)\right\}\right)}\le M\sup_{n\in\N}\|F_n\|_{L^2(\O)}^2,
\end{multline*}
for any $n\in\N.$ The same is true if we replace the conditions on $(a,b,c)$ by $(a,b,c)\in\A\times\A\times\{0\}$ satisfies~\eqref{ab} and
$\delta\le\delta_\star,$ where $\delta_\star$ is given by Lemma~$\ref{lemAB}.$ In this case, $M=M(|a|,|b|).$
\end{cor}

\begin{proof*}
Choosing $u_\ell^n$ and $\vi u_\ell^n$ as test functions, we obtain
\begin{multline}
\label{proofcorbound2-1}
\|\nabla u_\ell^n\|_{L^2(\O)}^2+\Re(a)\left(\|u_\ell^n\|_{L^{m+1}(\{|u_\ell^n|\le\ell\})}^{m+1}+\ell^m\|u_\ell^n\|_{L^1(\{|u_\ell^n|>\ell\})}\right)		\\
+\big(\Re(b)-\|V\|_{L^\infty(\O)}^2|\Re(c)|\big)\left(\|u_\ell^n\|_{L^2(\{|u_\ell^n|\le\ell\})}^2
+\ell\|u_\ell^n\|_{L^1(\{|u_\ell^n|>\ell\})}\right)\le\Re\int_{\O}F_n\ovl{u_\ell^n}\d x,
\end{multline}
\begin{multline}
\label{proofcorbound2-2}
\Im(a)\left(\|u_\ell^n\|_{L^{m+1}(\{|u_\ell^n|\le\ell\})}^{m+1}+\ell^m\|u_\ell^n\|_{L^1(\{|u_\ell^n|>\ell\})}\right)
		+\Im(b)\left(\|u_\ell^n\|_{L^2(\{|u_\ell^n|\le\ell\})}^2+\ell\|u_\ell^n\|_{L^1(\{|u_\ell^n|>\ell\})}\right)			\\
+\Im(c)\left(\|Vu\|_{L^2(\{|u_\ell^n|\le\ell\}))}^2+\ell\|V^2u\|_{L^1(\{|u_\ell^n|>\ell\}))}\right)=\Im\int_{\O}F_n\ovl{u_\ell^n}\d x,
\end{multline}
for any $(n,\ell)\in\N^2.$ If $(a,b,c)\in\A\times\A\times\{0\}$ satisfies \eqref{ab}, then we obtain
\begin{multline}
\label{proofcorbound2-3}
\|\nabla u_\ell^n\|_{L^2(\O)}^2+\delta\|u_\ell^n\|_{L^2(\O)}^2
+\Re(a)\left(\|u_\ell^n\|_{L^{m+1}(\{|u_\ell^n|\le\ell\})}^{m+1}+\ell^m\|u_\ell^n\|_{L^1(\{|u_\ell^n|>\ell\})}\right)		\\
+\big(\Re(b)-\delta\big)\left(\|u_\ell^n\|_{L^2(\{|u_\ell^n|\le\ell\})}^2+\ell\|u_\ell^n\|_{L^1(\{|u_\ell^n|>\ell\})}\right)=\Re\int_{\O}F_n\ovl{u_\ell^n}\d x,
\end{multline}
\begin{multline}
\label{proofcorbound2-4}
\Im(a)\left(\|u_\ell^n\|_{L^{m+1}(\{|u_\ell^n|\le\ell\})}^{m+1}+\ell^m\|u_\ell^n\|_{L^1(\{|u_\ell^n|>\ell\})}\right)					\\
+\Im(b)\left(\|u_\ell^n\|_{L^2(\{|u_\ell^n|\le\ell\})}^2+\ell\|u_\ell^n\|_{L^1(\{|u_\ell^n|>\ell\})}\right)=\Im\int_{\O}F_n\ovl{u_\ell^n}\d x,
\end{multline}
for any $(n,\ell)\in\N^2.$ For this last case, it follows from Lemma~\ref{lemAB}, Hölder's inequality and \eqref{young} that
\begin{multline*}
\|\nabla u_\ell^n\|_{L^2(\O)}^2+\frac{L}2\|u_\ell^n\|_{L^2(\{|u_\ell^n|\le\ell\})}^2
+L\big\|u_\ell^n\big\|_{L^{m+1}\left(\left\{\left|u_\ell^n\right|\le\ell\right\}\right)}^{m+1} \\
+\left(L\ell-M\|F_n\|_{L^\infty(\O)}\right)\|u_\ell^n\|_{L^1(\{|u_\ell^n|>\ell\})}
\le\frac{M^2}{2L}\|F_n\|_{L^2(\O)}^2.
\end{multline*}
Then the result follows by choosing for each $n\in\N,$ $\vphi(n)\in\N$ large enough to have $L\vphi(n)-M\|F_n\|_{L^\infty(\O)}\ge1.$ Now we turn out to the case \eqref{proofcorbound2-1}--\eqref{proofcorbound2-2}. Let $M$ and $A$ be given by Lemma~\ref{lemeap2} with $R=\|V\|_{L^\infty(\O)}.$ For each $n\in\N,$ let $\vphi(n)\in\N$ be large enough to have $\vphi(n)\ge2A\|F_n\|_{L^\infty(\O)}+1,$ if $|\o|<|\O|$ and $\vphi(n)=n,$ if $|\o|=|\O|.$ For each
$n\in\N,$ with help of \eqref{proofcorbound2-1} and \eqref{proofcorbound2-2}, we may apply Lemma~\ref{lemeap2} to $u_{\vphi(n)}^n$ with
$\o=\left\{x\in\O;\left|u_{\vphi(n)}^n(x)\right|\le\vphi(n)\right\},$ $\alpha=\vphi(n)^m,$ $\beta=\vphi(n)$ and $R=\|V\|_{L^\infty(\O)}.$ Hence the result.
\medskip
\end{proof*}

\section{Proofs of the main results}
\label{proofs}

\begin{vproof}{of Theorem~\ref{thmreg}.}
Property~1) follows from Proposition~4.5 in Bégout and D\'iaz~\cite{MR2876246} while Property~2) comes from Remark~4.7 in Bégout and
D\'iaz~\cite{MR2876246}. It remains to establish Property~3). Assume first that $F\in L^p(\O)$ and $V\in\bigcap\limits_{1<r<\infty} L^r(\O).$ It follows from the equation that for any $\eps\in(0,q-1),$ $\Delta u\in L^{q-\eps}(\O).$ We now recall an elliptic regularity result. If for some $1<s<\infty,$
$u\in L^s(\O)$ satisfies $\Delta u\in L^s(\O)$ and $\gamma(\nabla u.\nu)=0$ then $u\in W^{2,s}(\O)$ (Proposition~2.5.2.3, p.131, in Grisvard~\cite{MR3396210}). Since for any $\eps\in(0,q-1),$ $u,\Delta u\in L^{q-\eps}(\O)$ and $\gamma(\nabla u.\nu)=0$ (by assumption), by following the bootstrap method of the proof p.52 of Property~1) of Proposition~4.5 in Bégout and D\'iaz~\cite{MR2876246}, we obtain the result. Indeed, therein, it is sufficient to apply the global regularity result in Grisvard~\cite{MR3396210} (Proposition~2.5.2.3, p.131) in place of the local regularity result in Cazenave~\cite{caz-sle} (Proposition~4.1.2, p.101-102). Now, you turn out to the Hölder regularity. Assume $F\in C^{0,\alpha}(\ovl\O)$ and $V\in C^{0,\alpha}(\ovl\O).$ By global smoothness property in $W^{2,p}$ proved above, we know that $u\in W^{2,N+1}(\O)$ and $\gamma(\nabla u.\nu)=0$ in $L^{N+1}(\Gamma).$ It follows from the Sobolev's embedding, $W^{2,N+1}(\O)\inj C^{1,\frac1{N+1}}(\ovl\O)\inj C^{0,1}(\ovl\O),$ that for any $x\in\Gamma,$
$\frac{\partial u}{\partial\nu}(x)=0$ and $u\in C^{0,1}(\ovl\O).$ A straightforward calculation yields,
\begin{gather*}
\forall(x,y)\in\ovl\O^2, \; \left||u(x)|^{-(1-m)}u(x)-|u(y)|^{-(1-m)}u(y)\right|\le5|u(x)-u(y)|^m\le5|x-y|^m.
\end{gather*}
Setting, $g=F-(a|u|^{-(1-m)}u+(b-1)u+cVu),$ we deduce that $g\in C^{0,\alpha}(\ovl\O).$ Let $v\in C^{2,\alpha}(\ovl\O)$ be the unique solution to
\begin{gather*}
	\begin{cases}
			                  -\Delta v+v=g,		&	\text{ in } \O,	\\
		\dfrac{\partial v}{\partial\nu}=0,		&	\text{ on } \Gamma,
	\end{cases}
\end{gather*}
(see, for instance, Theorem~3.2 p.137 in Ladyzhenskaya and Ural'tseva~\cite{MR0244627}). It follows that $u$ and $v$ are two $H^1$-solutions of the above equations and since uniqueness holds in $H^1(\O)$ (Lax-Milgram's Theorem), we deduce that $u=v.$ Hence $u\in C^{2,\alpha}(\ovl\O).$ This concludes the proof\footnote{More directly, we could have said that since $u\in W^{2,N+1}(\O),$ $\gamma(\nabla u.\nu)=0$ and $\Delta u\in C^{0,\alpha}(\ovl\O)$ (by the estimate of the nonlinearity) then by Theorem~6.3.2.1, p.287, in Grisvard~\cite{MR3396210}, $u\in C^{2,\alpha}(\ovl\O).$ But this theorem requires $\O$ to have a $C^{2,1}$ boundary.}.
\medskip
\end{vproof}

\begin{vproof}{of Proposition~\ref{propreg}.}
We first establish Property~1). Since $\O$ has $C^{0,1}$ boundary and $u\in H^1_0(\O),$ it follows that $\gamma(u)=0.$ Moreover, Sobolev's embedding and equation \eqref{thmreg1} imply that $\Delta u\in L^2(\O).$ We then obtain that $u\in H^2(\O)$ (Grisvard~\cite{MR3396210}, Corollary~2.5.2.2,
p.131). Hence Property~1). We turn out to Property~2). It follows from equation~\eqref{thmreg1} that $\Delta u\in L^2(\O),$ so that \eqref{thmreg1} makes sense a.e. in $\O.$ Then Property~2) comes from the arguments of 2) of Remark~\ref{rmkdefsol}.
\medskip
\end{vproof}

\begin{lem}
\label{lemlaxmil}
Let $\vO\subset\R^N$ be a bounded open subset, let $V\in L^\infty(\O;\C),$ let $0<m<1,$ let $(a,b,c)\in\C^3$ and let $F\in L^2(\vO).$ Let
$\delta\in[0,1].$ Then for any $\ell\in\N,$ there exist a solution $u^1_\ell\in H^1_0(\vO)$ to
\begin{gather}
 \label{lemlaxmil1}
 -\Delta u_\ell+\delta u_\ell+f_\ell(u_\ell)=F, \text{ in } L^2(\vO),
\end{gather}
with boundary condition~\eqref{dir} and a solution $u^2_\ell\in H^1(\vO)$ to~\eqref{lemlaxmil1} with boundary condition~\eqref{neu} $($in this case, $\vO$ is assumed to have a $C^1$ boundary and $\delta>0),$ where
\begin{gather}
\label{lemlaxmil2}
\forall u\in L^2(\O), \; f_\ell(u)=
	\begin{cases}
		a|u|^{-(1-m)}u+(b-\delta)u+cV^2u,								&	\mbox{ if } |u|\le\ell, \medskip \\
		a\ell^m\dfrac{u}{|u|}+(b-\delta)\ell\dfrac{u}{|u|}+cV^2\ell\dfrac{u}{|u|},		&	\mbox{ if } |u| > \ell.
	\end{cases}
\end{gather}
If, in addition, $V$ is spherically symmetric then Symmetry Property~$\ref{sym}$ holds.
\end{lem}

\begin{proof*}
We proceed with the proof in two steps. Let $H=H^1_0(\vO),$ in the homogeneous Dirichlet case, and $H=H^1(\vO),$ in the homogeneous Neumann case. Let $\delta\in[0,1]$ \big(with additionally $\delta>0$ and $\Gamma$ of class $C^1$ if $H=H^1(\vO)\big).$ Step~1 below being obvious, we omit the proof. \\
{\bf Step~1.} $\forall G\in L^2(\vO),$ $\exists !u\in H$ s.t. $-\Delta u+\delta u=G.$ Moreover, $\exists\alpha>0$ s.t. $\forall G\in L^2(\vO),$
$\left\|(-\Delta+\delta I)^{-1}G\right\|_{H^1(\vO)}\le\alpha\|G\|_{L^2(\vO)}.$ Finally, Symmetry Property~\ref{sym} holds. \\
{\bf Step~2.} Conclusion. \\
For each $\ell\in\N,$ we define $g_\ell=-f_\ell+F\in C\big(L^2(\vO);L^2(\vO)\big).$ With help of the continuous and compact embedding
$i:H\inj L^2(\vO)$ and Step~1, we may define a continuous and compact sequence of mappings
$(T_\ell)_{\ell\in\N}$ of $H$ as follows. For any $\ell\in\N,$ set
\begin{gather*}
\begin{array}{rcccccl}
T_\ell:H		&     \stackrel{i}{\inj}	& 		             L^2(\vO)				& \xrightarrow{g_\ell}
			&          L^2(\vO)		&	\xrightarrow{(-\Delta+\delta I)^{-1}}		& H								\medskip \\
            u		&      \longmapsto	&                             i(u)=u					& \longmapsto
			&        g_\ell(u)      	&                        \longmapsto    				& (-\Delta+\delta u)^{-1}(g_\ell)(u)
\end{array}
\end{gather*}
Set $\rho=2\alpha(|a|+|b|+|c|+1)\left(\left(\|V\|_{L^\infty(\O)}^2+2\right)\ell|\vO|^\frac{1}{2}+\|F\|_{L^2(\vO)}\right).$ Let $u\in H.$ It follows that,
\begin{gather*}
\|T_\ell(u)\|_{H^1(\vO)}=\left\|(-\Delta+\delta I)^{-1}(g_\ell)(u)\right\|_{H^1(\vO)} \le\alpha\|g_\ell(u)\|_{L^2(\vO)}\le\rho.
\end{gather*}
Existence comes from the Schauder's fixed point Theorem applied to $T_\ell.$ The Symmetry Property~\ref{sym} is obtained by working in $H_\rad$ in place of $H$ \big(and in $H_\mathrm{even}$ and $H_\mathrm{odd}$ for $N=1\big).$
\medskip
\end{proof*}

\begin{vproof}{of Theorem~\ref{thmexista1}.}
Let  for any $u\in L^2(\O),$ $f(u)=a|u|^{-(1-m)}u+bu.$ Set $\O_n=\O\cap B(0,n).$ Let $(G_n)_{n\in\N}\subset\Dr(\O)$ be such that
$G_n\xrightarrow[n\to\infty]{L^2(\O)}F.$ Let $\big(u_\ell^n\big)_{(n,\ell)\in\N^2}\subset H^1_0(\O_n)$ a sequence of solutions of \eqref{lemlaxmil1} be given by Lemma~\ref{lemlaxmil} with $\vO=\O_n,$ $c=\delta=0$ and $F_n=G_{n|\O_n}.$ We define $\wt{u_\ell^n}\in H^1_0(\O)$ by extending $u_n$ by $0$ in $\O\cap\O_n^\c.$ We also denote by $\wt{f_\ell}$ the extension by $0$ of $f_\ell$ in $\O\cap\O_n^\c.$ By Corollary~\ref{corbound1}, there exists a diagonal extraction $\left(\wt{u_{\vphi(n)}^n}\right)_{n\in\N}$ of $\big(\wt{u_\ell^n}\big)_{(n,\ell)\in\N^2}$ which is bounded in $H^1_0(\O).$ By reflexivity of $H^1_0(\O),$ Rellich-Kondrachov's Theorem and converse of the dominated convergence theorem, there exist $u\in H^1_0(\O)$ and
$g\in L^2_\loc(\O;\R)$ such that, up to a subsequence that we still denote by $\left(\wt{u_{\vphi(n)}^n}\right)_{n\in\N},$
$\wt{u_{\vphi(n)}^n}\xrightarrow[n\to\infty]{L^2_\loc(\O)}u,$ $\wt{u_{\vphi(n)}^n}\xrightarrow[n\to\infty]{\text{a.e.~in }\O}u$ and
$\left|\wt{u_{\vphi(n)}^n}\right|\le g,$ a.e.~in $\O,$ By these two last estimates,
$\wt{f_{\vphi(n)}}\Big(\wt{u_{\vphi(n)}^n}\Big)\xrightarrow[n\to\infty]{\text{a.e.~in }\O}f(u)$ and
$\left|\wt{f_{\vphi(n)}}\Big(\wt{u_{\vphi(n)}^n}\Big)\right|\le C(g^m+g)\in L^2_\loc(\O),$ a.e.~in $\O.$ From the dominated convergence Theorem,
$\wt{f_{\vphi(n)}}\Big(\wt{u_{\vphi(n)}^n}\Big)\xrightarrow[n\to\infty]{L^2_\loc(\O)}f(u).$ Let $\vphi\in\Dr(\O).$ Let $n_\star\in\N$ be large enough to have
$\supp\vphi\subset\O_{n_\star}.$ We have by \eqref{lemlaxmil1},
\begin{gather*}
\forall n>n_\star, \;
\left\langle-\vi\Delta u_{\vphi(n)}^n+f_{\vphi(n)}\left(u_{\vphi(n)}^n\right)-F_n,\vphi_{|\O_n}\right\rangle_{\Dr^\p(\O_n),\Dr(\O_n)}=0.
\end{gather*}
The above convergencies lead to,
\begin{align*}
	&	\; \langle-\Delta u+f(u)-F,\vphi\rangle_{\Dr^\p(\O),\Dr(\O)}																\\
    =	&	\; \langle-u,\Delta\vphi\rangle_{\Dr^\p(\O),\Dr(\O)}+\langle f(u)-F,\vphi\rangle_{\Dr^\p(\O),\Dr(\O)}								\\
    =	&	\; \lim_{n\to\infty}\left\langle-\wt{u_{\vphi(n)}^n},\Delta\vphi\right\rangle_{\Dr^\p(\O),\Dr(\O)}
			+\lim_{n\to\infty}\left\langle\wt{f_{\vphi(n)}}\Big(\wt{u_{\vphi(n)}^n}\Big)-G_n,\vphi\right\rangle_{\Dr^\p(\O),\Dr(\O)}				\\
    =	&	\; \lim_{n\to\infty}\left\langle-\Delta u_{\vphi(n)}^n+f_{\vphi(n)}\left(u_{\vphi(n)}^n\right)
			-F_n,\vphi_{|\O_n}\right\rangle_{\Dr^\p(\O_n),\Dr(\O_n)}	\\
    =	&	\; 0.
\end{align*}
By density, we then obtain that $u\in H^1_0(\O)$ is a solution to $-\Delta u+f(u)=F, \text{ in } L^2(\O).$ Finally, if $F$ is spherically symmetric then $u$ (obtained as a limit of  solutions given by Lemma~\ref{lemlaxmil}) is also spherically symmetric. For $N=1,$ this includes the case where $F$ is an even function.
\medskip
\end{vproof}

\begin{vproof}{of Theorems~\ref{thmbound1} and \ref{thmbound3}.}
Choosing $u$ and $\vi u$ as test functions, we obtain
\begin{gather*}
	\|\nabla u\|_{L^2(\O)}^2+\Re(a)\|u\|_{L^{m+1}(\O)}^{m+1}+\Re(b)\|u\|_{L^2(\O)}^2=\Re\int_\O F\ovl u\d x,	\\
	\Im(a)\|u\|_{L^{m+1}(\O)}^{m+1}+\Im(b)\|u\|_{L^2(\O)}^2=\Im\int_\O F\ovl u\d x.
\end{gather*}
Theorem~\ref{thmbound1} follows immediately from Lemma~\ref{lemeap1} applied with $\o=\O,$ while Theorem~\ref{thmbound3} is a consequence of Lemma~\ref{lemAB} applied with $\delta=0$ and \eqref{young}. This ends the proof.
\medskip
\end{vproof}

\begin{vproof}{of Theorem~\ref{thmbound2}.}
Choosing $u$ and $\vi u$ as test functions, we obtain
\begin{gather*}
	\|\nabla u\|_{L^2(\O)}^2+\Re(a)\|u\|_{L^{m+1}(\O)}^{m+1}+\left(\Re(b)-|\Re(c)|\|V\|_{L^\infty(\O)}^2\right)\|u\|_{L^2(\O)}^2\le\int_\O|Fu|\d x,	\\
	|\Im(a)|\|u\|_{L^{m+1}(\O)}^{m+1}+|\Im(b)|\|u\|_{L^2(\O)}^2+|\Im(c)|\|Vu\|_{L^2(\O)}^2\le\int_\O|Fu|\d x.
\end{gather*}
The theorem follows Lemma~\ref{lemeap2} applied with $\o=\O,$ $R=\|V\|_{L^\infty(\O)}$ and $\alpha=\beta=0.$
\medskip
\end{vproof}

\begin{vproof}{of Theorems~\ref{thmexista2} and \ref{thmexista3}.}
We first assume that $\O$ is bounded. Let $H=H^1_0(\O),$ in the homogeneous Dirichlet case, and $H=H^1(\O),$ in the homogeneous Neumann case. Let $\delta_\star$ be given by Lemma~$\ref{lemAB}$ and let for any $u\in L^2(\O),$
$f(u)=a|u|^{-(1-m)}u+bu+cV^2u$ (with $c=0$ in the case of Theorem~\ref{thmexista3}). Let $(F_n)_{n\in\N}\subset\Dr(\O)$ be such that
$F_n\xrightarrow[n\to\infty]{L^2(\O)}F.$ Let $\big(u_\ell^n\big)_{(n,\ell)\in\N^2}\subset H$ a sequence of solutions of \eqref{lemlaxmil1} be given by Lemma~\ref{lemlaxmil} with $\vO=\O,$ $\delta=1$ for Theorem~\ref{thmexista2}, $\delta=\delta_\star$ for Theorem~\ref{thmexista3} and such $F_n.$ By Corollary~\ref{corbound2}, there exists a diagonal extraction $\left(u_{\vphi(n)}^n\right)_{n\in\N}$ of $\big(u_\ell^n\big)_{(n,\ell)\in\N^2}$ which is bounded in $W^{1,1}(\O)\cap\dot H^1(\O).$ Let $1<p<2$ be such that $W^{1,1}(\O)\inj L^p(\O).$ Then $\left(u_{\vphi(n)}^n\right)_{n\in\N}$ is bounded in $ W^{1,p}(\O)$ and there exist $u\in W^{1,p}(\O)\cap\dot H^1(\O)$ and $g\in L^p(\O;\R)$ such that, up to a subsequence that we still denote by
$\left(u_{\vphi(n)}^n\right)_{n\in\N},$ $u_{\vphi(n)}^n\xrightarrow[n\to\infty]{L^p(\O)}u,$
$\nabla u_{\vphi(n)}^n\weak\nabla u \text{ in } \left(L^2_\w(\O)\right)^N,$ as $n\tends\infty,$
$u_{\vphi(n)}^n\xrightarrow[n\to\infty]{\text{a.e.~in }\O}u,$ $\left|u_{\vphi(n)}^n\right|\le g,$ a.e.~in $\O$ and
$\left(u_{\vphi(n)}^n\1_{\left\{\left|u_{\vphi(n)}^n\right|\le\vphi(n)\right\}}\right)_{n\in\N}$ is bounded in $L^2(\O),$
where the last estimate comes from Corollary~\ref{corbound2}. By these three last estimates and Fatou's Lemma,
$u\in L^2(\O),$ $f_{\vphi(n)}\Big(u_{\vphi(n)}^n\Big)\xrightarrow[n\to\infty]{\text{a.e.~in }\O}f(u)-\delta u$ and
$\left|f_{\vphi(n)}\Big(u_{\vphi(n)}^n\Big)\right|\le C(g^m+g)\in L^p(\O),$ a.e.~in $\O.$ It follows that $u\in H^1(\O).$ From the dominated convergence Theorem, $f_{\vphi(n)}\Big(u_{\vphi(n)}^n\Big)\xrightarrow[n\to\infty]{L^p(\O)}f(u)-\delta u.$ Consider the Dirichlet boundary condition. We recall a Gagliardo-Nirenberg's inequality.
\begin{gather*}
\forall w\in H^1_0(\O), \; \|w\|_{L^2(\O)}^{N+2}\le C\|w\|_{L^1(\O)}^2\|\nabla w\|_{L^2(\O)}^N,
\end{gather*}
where $C=C(N).$ In particular, $C$ does not depend on $\O.$ Since $\left(u_{\vphi(n)}^n\right)_{n\in\N}\subset H^1_0(\O)$ is bounded in
$W^{1,1}(\O)\cap\dot H^1(\O),$ it follow from the above Gagliardo-Nirenberg's inequality that $\left(u_{\vphi(n)}^n\right)_{n\in\N}$ is bounded in $H^1_0(\O),$ so that $u\in H^1_0(\O).$ Now, we show that $u\in H$ is a solution. Let $m_0\in\N$ be large enough to have
$H^{m_0}(\O)\inj L^{p^\p}(\O).$ Let $v\in\Dr(\O),$ if $H=H^1_0(\O)$ and let $v\in H^{m_0}(\O),$ if $H=H^1(\O).$ By \eqref{lemlaxmil1}, we have for any $n\in\N,$
\begin{multline}
 \label{pthmexista2-6}
 \left\langle\nabla u_{\vphi(n)}^n,\nabla v\right\rangle_{L^2(\O),L^2(\O)}
 +\left\langle\delta u_{\vphi(n)}^n+f_{\vphi(n)}\left(u_{\vphi(n)}^n\right),v\right\rangle_{L^p(\O),L^{p^\p}(\O)}	\\
 -\langle F_n,v\rangle_{L^2(\O),L^2(\O)}=0.
\end{multline}
Above convergencies lead to allow us to pass in the limit in \eqref{pthmexista2-6} and by density of $\Dr(\O)$ in $H^1_0(\O)$ and density of
$H^{m_0}(\O)$ in $H^1(\O)$ (see, for instance, Corollary~9.8, p.277, in Brezis~\cite{MR2759829}), it follows that
\begin{gather*}
\forall v\in H, \;
\langle\nabla u,\nabla v\rangle_{L^2(\O),L^2(\O)}+\langle f(u),v\rangle_{L^2(\O),L^2(\O)}=\langle F,v\rangle_{L^2(\O),L^2(\O)}.
\end{gather*}
This finishes the proof of the existence for $\O$ bounded. Approximating $\O$ by an exhaustive sequence of bounded sets
$\left(\O\cap B(0,n)\right)_{n\in\N},$ the case $\O$ unbounded can be treated in the same way as in the proof of Theorem~\ref{thmexista1}. The symmetry property also follows as in the proof of Theorem~\ref{thmexista1}.
\medskip
\end{vproof}

\begin{vproof}{of Theorem~\ref{thmuni}.}
\label{proofthmuni}
Let $u_1,u_2\in H^1(\O)\cap L^{m+1}(\O)$ be two solutions of~\eqref{2} such that $Vu_1,Vu_2\in L^2(\O).$ We set $u=u_1-u_2,$ $f(v)=|v|^{-(1-m)}v$ and $g(v)=af(v)+bv+cV^2v.$ From Lemma~9.1 in Bégout and D\'iaz~\cite{MR2876246}, there exists a positive constant $C$ such that,
\begin{gather}
\label{proofthmuni1}
C\vint_\omega\frac{|u_1(x)-u_2(x)|^2}{(|u_1(x)|+|u_2(x)|)^{1-m}}\d x\le\langle f(u_1)-f(u_2),u_1-u_2\rangle_{L^\frac{m+1}{m}(\O),L^{m+1}(\O)},
\end{gather}
where $\omega=\Big\{x\in\O;|u_1(x)|+|u_2(x)|>0\Big\}.$ We have that $u$ satisfies $-\Delta u+g(u_1)-g(u_2)=0.$ Choosing $v=au$ as a test function, we get
\begin{gather*}
\Re(a)\|\nabla u\|_{L^2}^2+|a|^2\langle f(u_1)-f(u_2),u_1-u_2\rangle_{L^\frac{m+1}{m},L^{m+1}}+\Re(a\ovl b)\|u\|_{L^2}^2
+\Re\left(a\ovl c\right)\|Vu\|_{L^2}^2=0.
\end{gather*}
It follows from the above estimate and~\eqref{proofthmuni1} that,
\begin{gather*}
\Re(a)\|\nabla u\|_{L^2}^2+C|a|^2\vint_\omega\frac{|u_1(x)-u_2(x)|^2}{(|u_1(x)|+|u_2(x)|)^{1-m}}\d x+\Re(a\ovl b)\|u\|_{L^2}^2
+\Re\left(a\ovl c\right)\|Vu\|_{L^2}^2\le0,
\end{gather*}
which yields Property~1). Properties~2) and 3) follow in the same way.
\medskip
\end{vproof}

\begin{rmk}
\label{rmkuni}
It is not hard to adapt the above proof to find other criteria of uniqueness.
\end{rmk}

\section{On the existence of solutions of the Dirichlet problem for data beyond $\bs{L^2(\O)}$}
\label{diaz}

In this section we shall indicate how some of the precedent results of this paper can be extended to some data $F$ which are not in $L^2(\O)$ but in the more general Hilbert space $L^2(\O;\delta^\alpha),$ where $\delta(x)=\dist(x,\Gamma)$ and $\alpha\in(0,1).$
\medskip \\
In order to justify the associated notion of solution, we start by assuming that a function $u$ solves equation
\begin{gather}
\label{secdiaz}
-\Delta u+f(u)=F, \; \text{ in } \; \O,
\end{gather}
with the Dirichlet boundary condition~\eqref{dir}, $u_{|\Gamma}=0,$ and we multiply (formally) by $\ovl{v(x)}\delta(x),$ with $v\in H^1_0(\O;\delta^\alpha)$ \big(the weighted Sobolev space associated to the weight $\delta^\alpha(x)\big),$ we integrate by parts (by Green's formula) and we take the real part. Then we get,
\begin{gather}
\label{diaz1}
\Re\vint_\O\nabla u.\ovl{\nabla v}\,\delta^\alpha\d x+\Re\vint_\O\ovl v\,\nabla u.\nabla\delta^\alpha\d x	
+\Re\vint_\O f(u)\,\ovl v\,\delta^\alpha\d x=\Re\vint_\O F\,\ovl v\,\delta^\alpha\d x.
\end{gather}
To give a meaning to the condition~\eqref{diaz1}, we must assume that
\begin{gather}
\label{diaz2}
F\in L^2(\O;\delta^\alpha),
\end{gather}
where $\|F\|^2_{L^2(\O;\delta^\alpha)}=\dsp\vint_\O|F(x)|^2\delta^\alpha(x)\d x,$ and to include in the definition of solution the conditions
\begin{gather}
\label{diaz3}
u\in H^1_0(\O;\delta^\alpha)		\quad	\text{and}		\quad	f(u)\in L^2(\O;\delta^\alpha).
\end{gather}
The justification of the second term in~\eqref{diaz1} is far to be trivial and requires the use of a version of the following Hardy type inequality,
\begin{gather}
\label{diaz4}
\vint_\O|v(x)|^2\delta^{-(2-\alpha)}(x)\d x\le C\vint_\O|\nabla v(x)|^2\delta^\alpha(x)\d x,
\end{gather}
which holds for some constant $C$ independent of $v,$ for any $v\in H^1_0(\O;\delta^\alpha)$ once we assume that
\begin{gather}
\label{diaz5}
\O \text{ is a bounded open subset of } \R^N \text{ with Lipschitz boundary}
\end{gather}
(see, e.g., Kufner~\cite{MR802206} and also Dr\'abek, Kufner and Nicolosi~\cite{dkn}, Kufner and Opic~\cite{MR1069756}, Kufner and
Sänding~\cite{MR926688} and Ne\v{c}as~\cite{MR0163054}). Notice that under~\eqref{diaz5}, we know that $\delta\in W^{1,\infty}(\O)$ and so
\begin{gather*}
\left|\vint_\O\ovl v\,\nabla u.\nabla\delta^\alpha\d x\right|
=\left|\vint_\O\left(\delta^\frac\alpha2\nabla u\right).\left(\frac{\ovl v}{\delta^\frac\alpha2}\nabla\delta^\alpha\right)\d x\right|
\le\alpha\|\nabla\delta\|_{L^\infty(\O)}\|\nabla u\|_{L^2(\O;\delta^\alpha)}\|v\|_{L^2(\O;\delta^{-(2-\alpha)})}<\infty,
\end{gather*}
by Cauchy-Schwarz's inequality and~\eqref{diaz4}.
\begin{defi}
\label{defdiaz1}
Assumed~\eqref{diaz2},~\eqref{diaz5} and $\alpha\in(0,1),$ we say that $u\in H^1_0(\O;\delta^\alpha)$ is a \textit{solution} of \eqref{secdiaz} and \eqref{dir} in $H^1_0(\O;\delta^\alpha)$ if~\eqref{diaz3} holds and the integral condition~\eqref{diaz1} holds for any $v\in H^1_0(\O;\delta^\alpha).$
\end{defi}

\begin{rmk}
\label{rmkdiaz2}
Notice that $H^1_0(\O;\delta^\alpha)\inj L^2(\O)$ (by the Hardy's inequality~\eqref{diaz4} and \eqref{diaz5}). Moreover, since
\begin{gather}
\label{diaz6}
\delta^{-s\alpha}\in L^1(\O), \text{ for any } s\in(0,1),
\end{gather}
we know (Dr\'abek, Kufner and Nicolosi~\cite{dkn}, p.30) that
\begin{gather*}
H^1_0(\O;\delta^\alpha)\inj W^{1,p_s}(\O), \text{ with } \, p_s=\frac{2s}{s+1}.
\end{gather*}
\end{rmk}

\begin{rmk}
\label{rmkdiaz3}
Obviously, there are many functions $F$ such that $F\in L^2(\O;\delta^\alpha)\setminus L^2(\O)$ (for instance, if $F(x)\sim\frac{1}{\delta(x)^\beta},$ for some $\beta>0,$ then $F\in L^2(\O;\delta^\alpha),$ if $\beta<\frac{\alpha+1}{2}$ but $F\not\in L^2(\O),$ once $\beta\ge\frac12.$ This fact is crucial when the nonlinear term $f(u)$ involves a singular term of the form as in~\eqref{2} but with $m\in(-1,0)$ (see D\'iaz, Hern\'andez and
Rakotoson~\cite{MR2831448} for the real case).
\end{rmk}

\begin{rmk}
\label{rmkdiaz4}
We point out that in most of the papers dealing with weighted solutions of semilinear equations, the notion of solution is not justified in this way but merely by replacing the Laplace operator by a bilinear form which becomes coercive on the space $H^1_0(\O;\delta^\alpha).$ The second integral term in~\eqref{diaz1} is not mentioned \big(since, formally, the multiplication of the equation is merely by $v\in H^1_0(\O;\delta^\alpha)\big)$ but then it is quite complicated to justify that such alternative solutions satisfy the pde equation~\eqref{2} when they are assumed, additionally, that
$\Delta u\in L^2_\loc(\O).$ We also mention now (although it is a completely different approach) the notion of $L^1(\O;\delta)$-very weak solution developed recently for many scalars semilinear equations: see, e.g., Brezis, Cazenave, Martel and Ramiandrisoa~\cite{MR1357955}, D\'iaz and Rakotoson~\cite{MR2629573} and the references therein).
\end{rmk}

\noindent
By using exactly the same \textit{a priori} estimates, but now adapted to the space $H^1_0(\O;\delta^\alpha),$ we get the following result.

\begin{thm}
\label{thmdiaz}
Let $\O$ be a bounded open subset with Lipschitz boundary, $V\in L^\infty(\O;\R),$ $0<\alpha<1,$ $0<m<1,$ $(a,b,c)\in\C^3$ as in
Theorem~$\ref{thmexista2}$ and let $F\in L^2(\O;\delta^\alpha).$ Then we have the following result.
\begin{enumerate}
 \item[]
  \begin{enumerate}[$1)$]
   \item
    There exists at least a solution $u\in H^1_0(\O;\delta^\alpha)$ to \eqref{2}. Furthermore, any such solution belongs to $H^2_\loc(\O).$
   \item
    If, in addition, we assume the conditions of Theorem~$\ref{thmuni},$ this solution is unique in the class of $H^1_0(\O;\delta^\alpha)$-solutions.
  \end{enumerate}
\end{enumerate}
\end{thm}

\begin{rmk}
\label{rmkdiaz6}
In the proof of the \textit{a priori} estimates, it is useful to replace the weighted function $\delta$ by a more smooth function having the same behavior near $\Gamma.$ This is the case, for instance of the first eigenfunction $\vphi_1$ of the Laplace operator,
\begin{gather*}
\begin{cases}
-\Delta\vphi_1=\la_1\vphi_1,	&	\text{in } 	\O,		\\
             \vphi_{1|\Gamma}=0,	&	\text{on }	\Gamma.
\end{cases}
\end{gather*}
It is well-known that $\vphi_1\in W^{2,\infty}(\O)\cap W^{1,\infty}_0(\O)$ and that $C_1\delta(x)\le\vphi_1(x)\le C_2\delta(x),$ for any $x\in\O,$ for some positive constants $C_1$ and $C_2,$ independent of $x.$ Now, with this new weighted function, it is easy to see that the second term in~\eqref{diaz1} does not play any important role since, for instance, when taking $v=u$ as test function, we get that
\begin{align*}
	&	\; \Re\vint_\O\ovl u\,\nabla u.\nabla\vphi_1^\alpha\d x
			=\frac12\vint_\O\nabla|u|^2.\nabla\vphi_1^\alpha\d x=-\frac12\vint_\O|u|^2\Delta\vphi_1^\alpha\d x							\\
   =	&	\; \frac{\alpha\la_1}{2}\vint_\O|u|^2\vphi_1^\alpha\d x+\frac{\alpha(1-\alpha)}{2}\vint_\O|u|^2\vphi_1^{-(2-\alpha)}|\nabla\vphi_1|^2\d x\ge0.
\end{align*}
\end{rmk}

\section{Conclusions}
\label{conclusions}

In this section, we summarize the results obtained in Section~\ref{main} and give some applications.
\medskip \\
The next result comes from Theorems~\ref{thmexista1}, \ref{thmbound1} and \ref{thmuni}.

\begin{thm}
\label{thmnls1}
Let $\O$ an open subset of $\R^N$ be such that $|\O|<\infty$ and assume $0<m<1,$ $(a,b)\in\C^2$ and $F\in L^2(\O).$ Assume that $\Re(b)>-\frac1{\CP^2}$or $\Im(b)\neq0,$ where $\CP$ is the Poincaré's constant in~\eqref{poincare}. Then there exists at least a solution $u\in H^1_0(\O)$ to
\begin{gather}
\label{eq1}
-\Delta u+a|u|^{-(1-m)}u+bu=F, \text{ in } L^2(\O).
\end{gather}
Furthermore, $\|u\|_{H^1_0(\O)}\le C(\|F\|_{L^2(\O)},|\O|,|a|,|b|,N,m).$ Finally, if $\Re(a)\ge0$ and $\vecteur a.\vecteur b\ge0$ then the solution is unique.
\end{thm}

\noindent
In the above theorem, the complex numbers $a$ and $b$ are seen as vectors $\vecteur a$ and $\vecteur b$ of $\R^2.$ Consequently, $\vecteur a.\vecteur b$ denotes the scalar product between these vectors of $\R^2.$

\medskip
\noindent
The novelty of Theorem~\ref{thmnls1} is about the range of $(a,b):$ we obtain existence of solution with, for instance, $(a,b)\in\R_-\times(-\eps,0),$ with $\eps>0$ small enough, or $(a,b)=(-1+\vi,-1-\vi).$ Recall that, up to today, existence was an open question when $(a,b)\in\R_-\times\R_-$ or
$[a,b]\cap\R_-\times\vi\{0\}\neq\emptyset$ (Bégout and D\'iaz~\cite{MR2876246}). Knowing that for such $(a,b)$ equation \eqref{eq1} admits solutions, it would be interesting if, whether or not, solutions with compact support exist, as in Bégout and D\'iaz~\cite{MR2876246}.

\medskip
\noindent
By Theorems~\ref{thmexista2}, \ref{thmbound2} and \ref{thmuni}, we get the following result.

\begin{thm}
\label{thmnls2}
Let $\O\subseteq\R^N$ be a bounded open subset, let $0<m<1$ and let $(a,b,c)\in\C^3$ be such that $\Im(a)<0,$ $\Im(b)<0$ and $\Im(c)\le0.$ For any $F\in L^2(\O),$ there exists at least a solution $u\in H^1(\O)$ to
\begin{gather}
\label{eq2}
-\Delta u+a|u|^{-(1-m)}u+bu+c|x|^2u=F, \text{ in } L^2(\O),
\end{gather}
with boundary condition \eqref{dir} or \eqref{neu}\textsuperscript{\textnormal{\ref{fn1}}}. Furthermore,
\begin{gather*}
\|u\|_{H^1(\O)}\le C(|a|,|b|,|c|)(R^2+1)\|F\|_{L^2(\O)},
\end{gather*}
where $B(0,R)\supset\O.$ Finally, if $\Re(a)\ge0,$ $\vecteur a.\vecteur b\ge0$ and $\vecteur a.\vecteur c\ge0$ then the solution is unique.
\end{thm}

\noindent
Since, now, we are able to show that equation \eqref{eq2} admits solutions, we can study the propagation support phenomena. Indeed, we can show that, under some suitable conditions, there exists a self-similar solution $u$ to
\begin{gather*}
\vi u_t+\Delta u=a|u|^{-(1-m)}u+f(t,x), \text{ in } \R^N,
\end{gather*}
such that for any $t>0,$ $\supp u(t)$ is compact (see~Bégout and D\'iaz~\cite{MR3193996}).

\medskip
\noindent
Now, we turn out to equation~\eqref{eq1} by extending some results found in Bégout and D\'iaz~\cite{MR2876246}. These results are due to Theorems~\ref{thmexista3}, \ref{thmbound3} and \ref{thmuni}.

\begin{thm}
\label{thmnls3}
Let $\O\subseteq\R^N$ be an open subset of $\R^N,$ let $0<m<1$ and let $(a,b)\in\A^2$ satisfies~\eqref{ab}. For any $F\in L^2(\O),$ there exists at least a solution $u\in H^1(\O)\cap L^{m+1}(\O)$ to
\begin{gather}
\label{eq3}
-\Delta u+a|u|^{-(1-m)}u+bu=F, \text{ in } L^2(\O)+L^\frac{m+1}{m}(\O),
\end{gather}
with boundary condition \eqref{dir} or \eqref{neu}\textsuperscript{\textnormal{\ref{fn1}}} $($in this last case, $\O$ is assumed bounded$)$. Furthermore,
\begin{gather*}
\|u\|_{H^1(\O)}^2+\|u\|_{L^{m+1}(\O)}^{m+1}\le M(|a|,|b|)\|F\|_{L^2(\O)}^2.
\end{gather*}
Finally, if $\Re(a)\ge0$ and $\vecteur a.\vecteur b\ge0$ then the solution is unique.
\end{thm}

\noindent
When $|\O|<\infty,$ Theorem~\ref{thmnls3} is an improvement of  Theorem~4.1 of Bégout and D\'iaz~\cite{MR2876246}, since we may choose
$F\in L^2(\O),$ instead of $F\in L^\frac{m+1}{m}(\O)$ and that $L^\frac{m+1}{m}(\O)\subsetneq L^2(\O).$ In addition, this existence result extends to the homogeneous Neumann boundary condition. In this context, we may show three kinds of new results, under assumptions of Theorem~\ref{thmnls3}.

$\bullet$ If $\O=\R^N$ and if $F\in L^2(\R^N)$ has compact support then equation~\eqref{eq3} admits solutions and any solution is compactly supported.

$\bullet$ If $\|F\|_{L^2(\O)}$ is small enough and if $F$ has compact support then equation~\eqref{eq3} admits solutions with the homogeneous Dirichlet boundary condition and any solution is compactly supported in $\O.$

$\bullet$ If $\|F\|_{L^2(\O)}$ is small enough, if $\vecteur a.\vecteur b>0$ and if $F$ has compact support then equation~\eqref{eq3} admits a unique solution with the homogeneous Neumann boundary condition and, in fact, this solution is compactly supported in $\O.$

\noindent
For more details, see~Bégout and D\'iaz~\cite{MR3190983}. Finally, in Section~\ref{diaz} we extended our techniques of proofs to the case in which the datum $F$ is very singular near the boundary of $\O$ but still is in some weighted Lebesgue space (see Theorem~\ref{thmdiaz}).

\medskip
\noindent
\textbf{Acknowledgements} \\
\baselineskip .5cm
The first author is very grateful to Professor Thierry Cazenave for helpful comments about the presentation of this paper and both authors thank the anonymous referee for stimulating remarks.

\baselineskip .4cm

\addcontentsline{toc}{section}{References}

\def\cprime{$'$}

\end{document}